\documentclass[11pt,a4paper]{article}
\usepackage[latin1]{inputenc}
\usepackage{indentfirst}
\usepackage{amsmath}
\usepackage{amsfonts}
\usepackage{amsthm}
\usepackage{amssymb}
\usepackage{graphics}
\usepackage{graphicx}

\def\H{\mathcal{H}}

\def\N{\mathbb{N}}

\def\P{\mathbb{P}}
\def\E{\mathbb{E}}

\def\R{\mathcal{R}}
\def\Rp{\overset{\rightarrow }{\scriptstyle{ \mathcal{R}}}}
\def\Rpt{\overset{\longrightarrow }{\scriptstyle{ \mathcal{R}(t)}}}

\def\Rmt{\overset{\longleftarrow }{\scriptstyle{ \mathcal{R}(t)}}}

\def\C{\mathcal{C}}
\def\RRR{\mathbb{R}}

\def\t{\textrm}
\def\w{\widetilde}

\def\B{\B }
\def\taup{\overset{\rightarrow }{\scriptstyle{\tau}}}
\def\taum{\overset{\leftarrow}{\scriptstyle{\tau}}}
\def\B{\mathbf{B} _0}
\def\ind{{\mathchoice {\rm 1\mskip-4mu l} {\rm 1\mskip-4mu l}
{\rm 1\mskip-4.5mu l} {\rm 1\mskip-5mu l}}}

\newcommand{\be} {\begin{equation}}
\newcommand{\ee} {\end{equation}}
\newcommand{\bea} {\begin{eqnarray}}
\newcommand{\eea} {\end{eqnarray}}
\newcommand{\Bea} {\begin{eqnarray*}}
\newcommand{\Eea} {\end{eqnarray*}}

\newtheorem{Thm}{Theorem}
\newtheorem{Lem}{Lemma}

\newtheorem{Pte}{Proposition}
\newtheorem{Cor}{Corollary}
\theoremstyle{definition} 
\theoremstyle{definition} \newtheorem*{key}{Key words}
\theoremstyle{definition} \newtheorem*{ams}{A.M.S. Classification}

\theoremstyle{remark}\newtheorem{Rque}{Remark}
\theoremstyle{remark}\newtheorem{fig}{Figure}
\theoremstyle{remark}\newtheorem{ex}{Example}

\begin{document}
\title{On a model for the storage of files on a hardware II : \\ Evolution of a typical data block.}
\author{Vincent Bansaye \footnote{Laboratoire de Probabilités et Modèles Aléatoires. Université Pierre et Marie Curie et C.N.R.S. UMR 7599. 175, rue du
Chevaleret, 75 013 Paris, France.
$\newline$
\emph{e-mail} : bansaye@ccr.jussieu.fr }}
\maketitle
\vspace{3cm}
\begin{abstract}We consider  the generalized version in continuous time of the parking problem
of Knuth introduced in \cite{vb}. Files arrive following a Poisson point process and are stored on a hardware
identified with the real line, at the right of their arrival point.
We study here the evolution of the extremities of the data block straddling $0$, which is empty at time $0$ and is equal to
 $\RRR$ at a deterministic time.
\end{abstract}
\begin{key}Parking problem. Data storage. Random covering.  Poisson point process. Lévy process.
\end{key}
\begin{ams}60D05, 60G55, 60J80, 68B15.
\end{ams}

\section{Introduction}

This paper is a continuation of $\cite{vb}$ but it can
be read independently. We consider  a generalized
 version in continuous time of the original parking problem of Knuth, as a model for the storage of files on a hardware. We
are interested in the evolution of a typical data block while files are stored on the hardware and we shall
characterize the process of the extremities and the length of this block.\\

We recall now the process of storage of files. In the original problem of Knuth, files arrive
successively  at location chosen
uniformly among  $n$ spots. They are stored in the first free spot
at the right of their arrival point (see \cite{chass, Flaj, Foa}).
In the model considered here, the hardware is  identified  with the real line and a file labelled
$i$ of length (or size) $l_i$ arrives at time $t_i$ on the real line
at location $x_i$.  The storage of this file uses the free portion of size $l_i$ of the real
line at the right of $x_i$ as close to $x_i$ as possible (see Figure \ref{fig1}). That is :
it covers $[x_i,x_i+l_i[$ if this interval  is free at time $t_i$.
Otherwise it is shifted to the right until a free space is found
and it may  be split into several parts which are stored in the closest
free spots.
$\newline$

\begin{fig}
\label{fig1}
Arrival and storage of a file on the hardware, where the data blocks are represented by black rectangles.$\qquad \qquad 
\qquad \qquad \qquad  \qquad \qquad \qquad  \qquad \qquad \qquad \qquad$
$\includegraphics[scale=0.43]{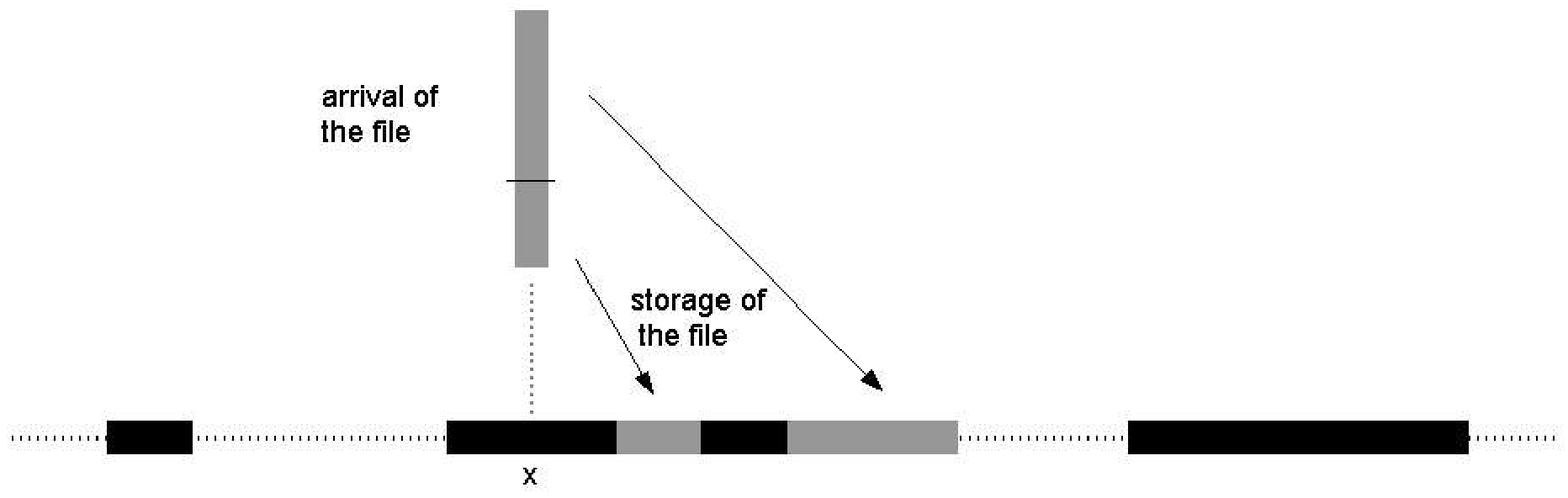}$
\end{fig}
$\newline$
The arrival of files follow  a Poisson point
process (PPP) :
 $\{(t_i,x_i,l_i)\  : \ i\in \N\}$ is a  PPP  with intensity $\t{d}t\otimes \t{d}x\otimes \nu(\t{d}l)$ on
$\RRR^+\times \RRR\times \RRR^+$. We  denote $\bar{\nu}(x)=\nu(]x,\infty])$ and we assume  
$\t{m}:=\int_0^\infty l \nu(\t{d}l)<\infty$. So
$\t{m}$ is the mean of the total sizes of files which arrive during a unit interval time on some interval
with unit length.
In \cite{vb}, this random covering has been constructed rigorously and  some statistics of this covering
were given. We proved that
the hardware becomes full at a deterministic time equal to $1/\t{m}$, studied the asymptotics at this saturation
time and  characterized the distribution of the covering at a fixed time  by  giving the joint
distribution of the block of data straddling $0$ and the free spaces on the sides of this block. \\ \\

In this work, we focus on the dynamics of the covering and we shall study the block of data straddling
a typical point, say $0$ for simplicity, which is denoted by $\B$. Thus $\B (t)$ is the block of data of the hardware
 containing $0$ at time $t$. We will show that its extremities and its length are pure jump Markov processes.\\
Specifically, if a   file arrives at time $t$ at the left of $\B (t-)$ and  cannot be stored
entirely at its left, it yields a jump of the left extremity of $\B$. The
data of this file  which cannot be stored at the left of $\B (t-)$ are called
\emph{remaining data}. These remaining data yield a jump of the right extremity of $\B$
(see Figure \ref{figg}). We shall prove that these events  happen
at instants which accumulate at  $1/\t{m}$ and induce a random partition of the time interval $[0,\t{1/m}]$ with the
Poisson-Dirichlet distribution (Theorem \ref{g}) and that the  jumps of the extremities at these instants
form  a PPP on $[0,1/\t{m}]\times \RRR_+\times \RRR_+$ (Proposition $\ref{pppd}$). Moreover  the successive
quantities of remaining data form an iid sequence (Corollary \ref{R}). \\
If a file arrives on $\B$, it yields a jump of the right
extremity  only (see Figure \ref{figgg}). The other files do not induce immediately a jump of $\B$ and
 we get the evolution of $(\B (t))_{t\geq 0}$ (Theorem
$\ref{couple}$). Finally, we prove that the process describing the length of $(\B (t))_{t\geq 0}$ is
a branching process with immigration (Corollary \ref{branch}). \\ \\

\begin{fig}
\label{figg}
Jumps of the extremities of $\B $ ($\Delta g (t)$ and $\Delta d (t)$) and remaining data induced by the arrival of a file at
 time $t$  at the left of $\B (t-) $. $\qquad \qquad \qquad \qquad \qquad$
$\includegraphics[scale=0.45]{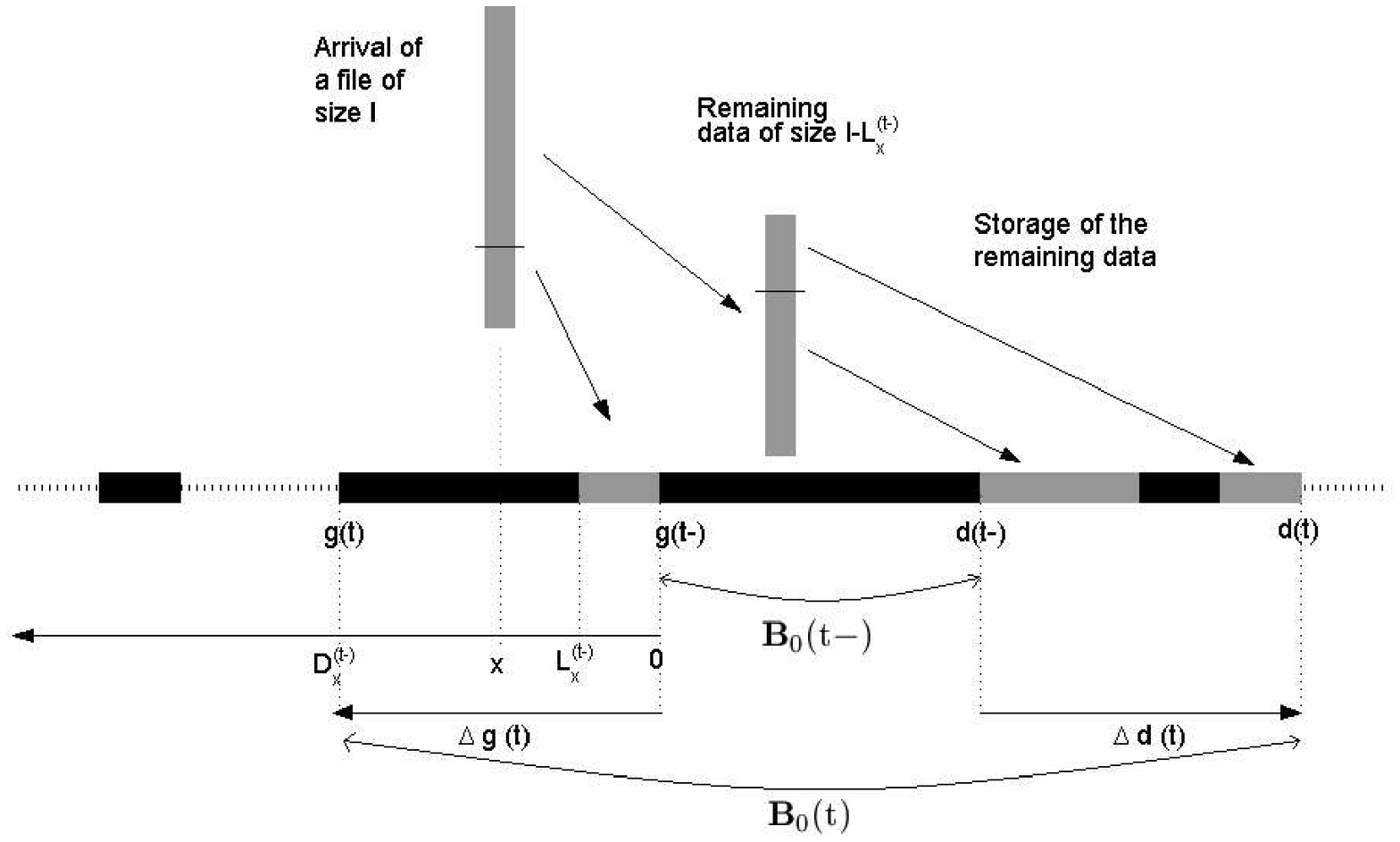}$
\end{fig}

\begin{fig}
\label{figgg}
Jump of the right extremity of  $\B $ ($\Delta d (t)$)  induced by the
arrival of a file at time $t$   on $\B (t-)$.
$\qquad \qquad \qquad \qquad \qquad \qquad \qquad \qquad \qquad \qquad \qquad \qquad \qquad$
$\includegraphics[scale=0.45]{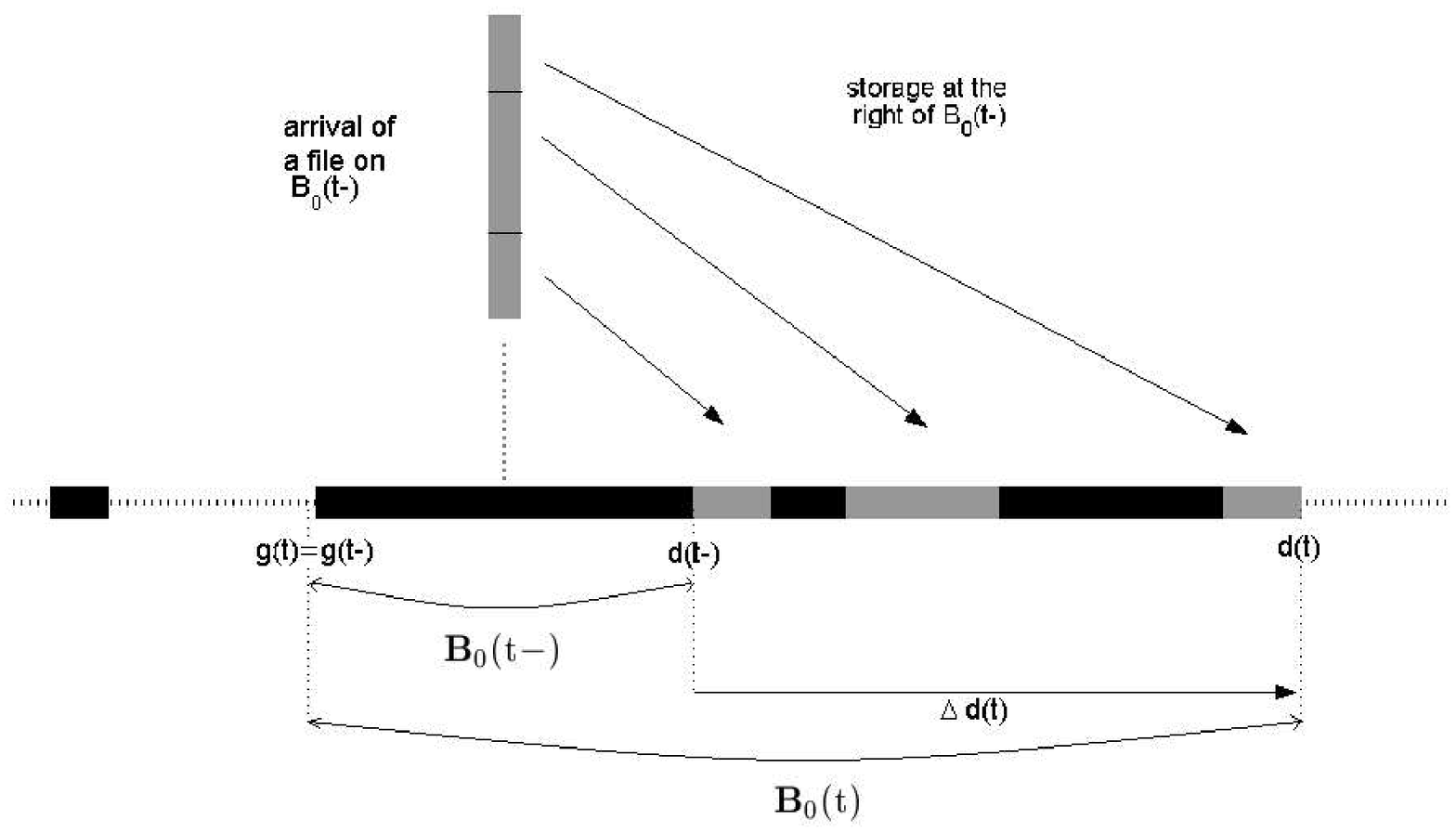}$
\end{fig}

\section{Preliminaries}

The covering $\C(t)$ described in Introduction has been
constructed in Section 2.1 in \cite{vb} and we recall some useful
results of this work. We denote by $\R(t)$ the complementary set
of $\C(t)$. It is natural and convenient to decide that files and so $\C(t)$ and $\R(t)$ are closed at the left, open at the
right. We  introduce the  process $(Y^{(t)}_x)_{x\in\RRR}$
defined by
\be
\label{defy}
Y^{(t)}_0:=0 \quad ; \quad  Y^{(t)}_b-Y^{(t)}_a= \sum_{\substack{t_i\leq t \\ x_i \in ]a,b]}} l_{i} \ \ -  \ \ (b-a) \quad \t{for} \ a<b.
\ee
It has càdlàg paths and stationary independent increments. The
process $(Y^{(t)}_x)_{x\geq 0}$ is then a Lévy process. Its
 drift  is equal to $-1$ and  its  Lévy measure is equal to $t\nu$. Its Laplace exponent $\Psi^{(t)}$ defined by
 \be
 \label{defexplap}
 \forall \ \rho\geq 0, \quad  \E(\t{exp}(-\rho Y^{(t)}_x))=\t{exp}(-x \Psi ^{(t)}(\rho )),
\ee
is given by
\be
\label{Psi}
\forall \  \rho\geq 0, \quad \Psi^{(t)}(\rho)=-\rho+\int_0^{\infty} \big( 1-e^{-\rho x} \big) t \nu(\t{d}x).
\ee
Introducing also its infimum process $I^{(t)}_x:=\t{inf}\{Y^{(t)}_y : y\leq x\}$ for every $x\in\RRR$, we got
the following  expression for the covering  and the free space
 \be \label{comp}
\C (t)= \{ x\in \RRR : \ Y^{(t)}_x > I^{(t)}_x \}, \qquad
\R(t)=\{x \in \RRR : Y^{(t)}_x=I^{(t)}_x\}\quad \t{a.s.} \ee
$\newline$
\begin{fig}Representation of $Y$ on a part of the hardware.
\begin{center}
$\includegraphics[scale=0.33]{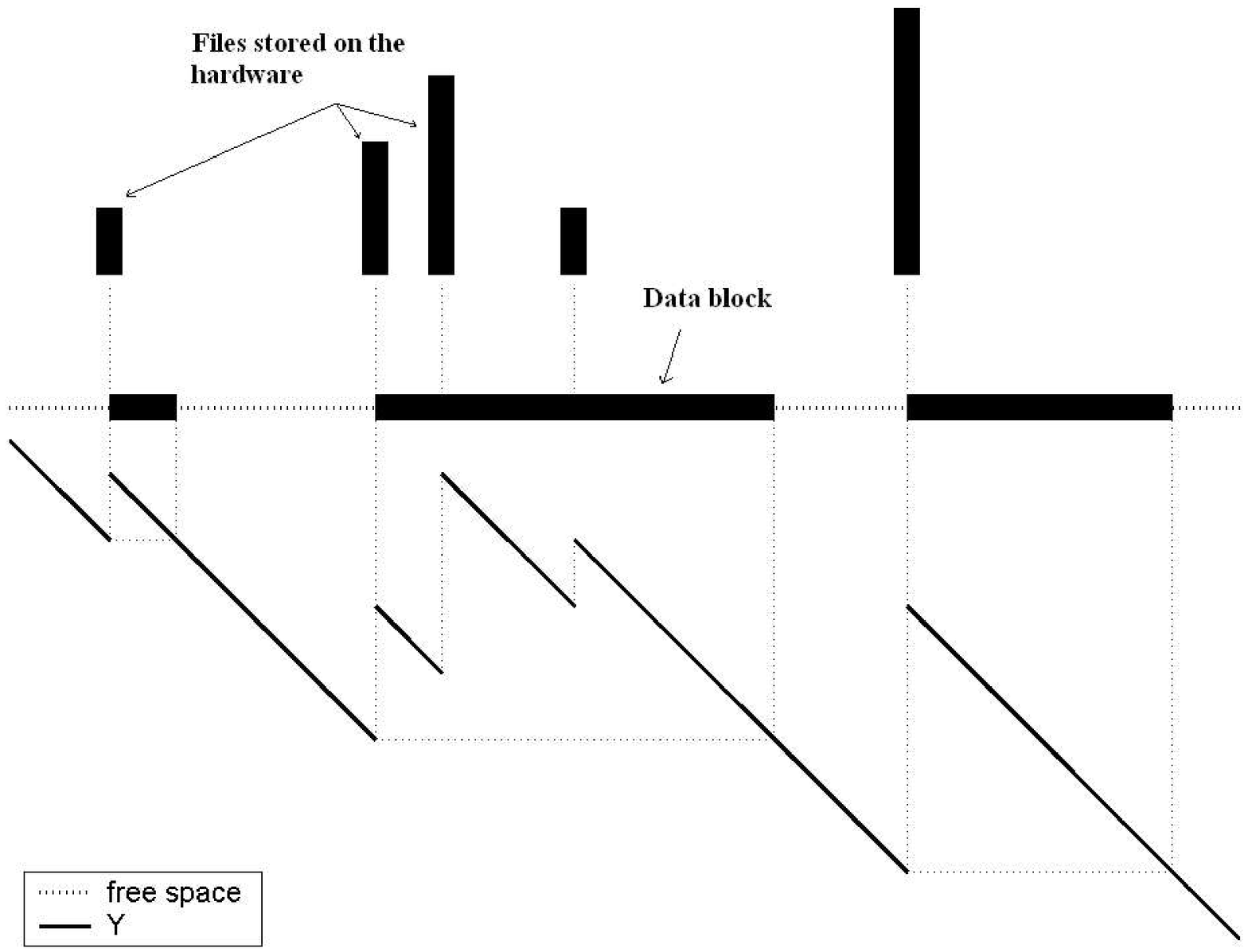}$
\end{center}
\end{fig}
$\newline$
The time when the hardware becomes full is equal to $1/\t{m}$, that is
a.s \ $\C(t)=\RRR$ iff $t\geq 1/$m. Thus we already know that $\B
(0)=\varnothing$ and $\B (1/\t{m})=\RRR$ and  we shall study $(\B
(t))_{t\in[0,1/\t{m}]}$. In that view, we introduce $g(t)$ (resp.
$d(t)$, resp. $l(t)$) the left extremity (resp. the right extremity,
resp. the length) of the data block containing $0$ :
$$\B (t)=[g(t),d(t)[, \qquad  l(t)=d(t)-g(t).$$
We will  also need  the free space at the  right of $\B(t)$ denoted by $\Rpt$ and  at the left of  $\B(t)$, turned
over, closed at the left and open at the right,  denoted by  $\Rmt$. If $\R\subset \RRR$ and $\R=\sqcup_{n\in\N}[a_n,b_n[$,
we denote by $\w{\R}=\sqcup_{n\in\N}[-b_n,-a_n[$ the symmetric set closed at the left and open at the right. Then we can
define (see Section 3 in \cite{vb}  for details)
\be
\Rpt :=(\R(t)-d(t))\cap [0,\infty], \qquad \Rmt := \overset{\longrightarrow }{\w{\scriptstyle{\mathcal{R}(t)}}}, \nonumber
\ee
which satisfy the following identity
\be
\label{reprR}
\R(t)=(d(t)+\Rpt)\sqcup(-\w{g(t)+\Rmt}).
\ee
In \cite{vb} Section 3, we proved that $\Rpt$ and  $\Rmt$ are   the range  of the processes $(\taup^{(t)}_x)_{x\geq 0}$ and $(\taum^{(t)}_x)_{x\geq 0}$ respectively
 defined by
$$\taup^{(t)}_x:=\t{inf}\{y\geq 0 : \vert\Rpt\cap [0,y]\vert>x\}, \qquad \taum^{(t)}_x:=\t{inf}\{y\geq 0 : \vert\Rmt\cap [0,y]\vert>x\}.$$
Moreover denoting by $\kappa^{(t)}$ the inverse function of $-\Psi^{(t)}$ and by $\Pi^{(t)}$ its Lévy measure :
\be
\label{kappat}
\kappa^{(t)}\circ (-\Psi^{(t)})=\rm{Id}, \qquad \forall \rho\geq 0, \quad 
\kappa ^{(t)}(\rho)= \rho  +\int_0^{\infty} (1-e^{-\rho x}) \Pi^{(t)} (\t{d}x),
\ee
enabled  us to describe  $\R(t)$ in the following way  :
$\newline$
\begin{Thm}
\label{loi}
(i) The processes $\taup^{(t)}$ and $\taum^{(t)}$ are two indepedent subordinators with 
Laplace exponent $\kappa^{(t)}$, which are independent of  $(g(t),d(t))$.\\
(ii) The  distribution of $(g(t),d(t))$ is specified by :
$$(g(t),d(t))=(-Ul(t),(1-U)l(t)),$$
$$\P(l(t) \in \emph{d}x)=(1-\emph{m}t)\big(\delta _0 (\emph{d}x)+\ind _{\{x>0\}} x\Pi^{(t)}(\emph{d}x) \big)$$
where $U$ uniform random variable on $[0,1]$ independent of $l(t)$.
\end{Thm}

$\newline$

For the basic example $\nu=\delta_1$, we got  for all $x\in \RRR_+$ and $n
\in \N$, 
\be
\label{diracY}
\P(Y_x^{(t)}+x=n)=e^{-tx}\frac{(tx)^n}{n!},
\ee
\be\label{dirac}
\P(\taup^{(t)}_x=x+n)=\frac{x}{x+n}e^{-t(x+n)}\frac{(t(n+x))^{n}}{n!}, \qquad
\Pi^{(t)}(n)= \frac{(tn)^n}{n .n!} e^{-tn}
\ee
Thus $l(t)$ follows a
size biased  Borel law :
$$\P(l(t)=n)=(1-t) \frac{(tn)^n}{n!} e^{-tn}.$$

$\newline$

We proved also the following identities :
\be
\label{rel}
\bar{\Pi}^{(t)}(0)=t\bar{\nu}(0), \qquad  \int_0^{\infty} x \Pi^{(t)}(\t{d}x)  =\frac{\t{m}t}{1-\t{m}t}, \qquad [\kappa^{(t)}]'(0)=\frac{1}{1-\t{m}t},
\ee
and the following identities of measures on $\RRR_+\times \RRR_+$,
\be
\label{egal}
x\P(\taum^{(t)}_l \in \t{d}x)\t{d}l=x\P(\taup^{(t)}_l \in \t{d}x)\t{d}l=l\P(-Y^{(t)}_x \in \t{d}l)\t{d}x.
\ee
Finally, we recall   a useful expression for the law of $g(t)$.  For all
$t\in [0,\t{1/m}[$ and  $\lambda \geq 0$,
\be
\label{gg}
\E\big(\t{exp}\big(\lambda g(t) \big) \big)=\t{exp}\bigg(\int_0^{\infty} (e^{-\lambda x}-1) x^{-1} \P(Y^{(t)}_x> 0)\t{d}x\bigg).
\ee
$\newline$

We can focus now  on the evolution of the block containing $0$, $\B$.
First, we prove some properties of absence of memory (Section 3)  : the
evolution of $\B$ after time $t$ depends from the past of this
block only through $l(t)$ (Markov property). Then we focus on the
left extremity  : it is an additive process and we give its Lévy
measure. As a consequence, we get the distribution of the instants at
which the left extremity jumps (Section 4). We then derive the distribution of the remaining data which
completes the description of the process of storage at the
left extremity (Section 5).  By taking also into account the data fallen on
$\B$, we get  then the evolution of $(g(t),d(t))$ (Section 6).
The latter characterizes the evolution of the right extremity and the
length (Section 7). $\newline$

\section{Markov property of $\B$}

We have already proved that $\R(t)$ enjoys a 'spatial' regeneration property (see Proposition 3 in \cite{vb}). To study the
evolution of $\B$, we need 'time' regeneration property. Here we
prove that
the evolution of the block containing $0$ up to time $t$ is independent
of the covering outside $[g(t),d(t)]$ up to time $t$. In Section 5, this property will ensure that
the evolution of $\B$ after time $t$ depends from the past of this block only through $l(t)$ (Markov property).

$\newline$
\begin{Pte}
\label{indd}
For every  $t \in [0,\emph{1/m}[$, the following three processes with values in the space of subsets of $\RRR$ \\
$. \qquad \qquad \qquad (g(t)-\R(s))\cap [0,\infty[, \quad \ \  0\leq s\leq t$, \\
$. \qquad \qquad \qquad  (\R(s)-d(t))\cap [0,\infty[, \quad \ \   0\leq s \leq t$, \\
$. \qquad \qquad \qquad  \R(s)\cap [g(t),d(t)], \quad  \quad \  \ \ \ 0\leq s\leq t,$ \\
 are independent.
\end{Pte}
\begin{Rque}
Actually, we have the following regeneration  property : $ \forall t \in [0,\t{1/m}[$, $\forall x \in \RRR$,   $\big((\R(s)-d_x(\R(t)))\cap [0,\infty[ : s \in [0,t]\big)$ is independent of
$\big((\R(s)-d_x(\R (t)))\cap ]-\infty,0] : s \in [0,t]\big)$
and is distributed as  $\big((\R(s)-d_0(\R(t)))\cap [0,\infty[ : s \in [0,t]\big)$.
\end{Rque}
$\newline$

This result is a direct consequence of the following lemma where we consider the point  processes of files  until time $t$
at the left of/at the right of/inside $[g,d]$ :
$$ P_g(t):=\{(t_i,g-x_i,l_i) :   t_i\leq t, \ x_i<g\},  \ \ \ \  P^d(t):=\{(t_i,x_i-d,l_i) :   t_i\leq t, \ d< x_i\}, $$
$$ P_g^d (t):=\{(t_i,x_i,l_i) :   t_i\leq t, \ g\leq x_i\leq d\}.  $$
\begin{Lem}
\label{indpp} For every $t \in [0,\emph{1/m}[$, the point processes $P_{g(t)}(t)$, $P_{g(t)}^{d(t)} (t)$ and $P_{d(t)} (t)$
 are independent.
\end{Lem}
\begin{proof} First we prove a weaker result, where times $(t_i)_{i\in\N}$ are not taken into account. Denote by
$(\w{Y}^{(t)}_x)_{x\geq 0}$  the càdlàg version
of $(Y^{(t)}_{-x})_{x\geq 0}$. This is a spectrally negative Lévy process  with bounded variation, which drifts to $\infty$. Note that,
 \bea
 g(t)&=&g_0(\R(t))= \t{sup}\{x\leq 0 : \ Y^{(t)}_x=I_x^{(t)}\} \nonumber \\
&=&\t{sup}\{x\leq 0 : \ Y^{(t)}_{x^-}=I^{(t)}_0\}=-\t{inf}\{x\geq 0 : \
\w{Y}^{(t)}_x=\t{inf}\{\w{Y}^{(t)}_z : z\geq 0\}\}. \nonumber \eea Then
$(\w{Y}^{(t)}_{-g(t)+x}-\w{Y}^{(t)}_{-g(t)})_{x\geq 0}$ is
independent of $(\w{Y}^{(t)}_x)_{0\leq x\leq -g(t)}$ (decomposition
of a Lévy process at its infimum \cite{mil}). Considering the locations and sizes of the jumps of these
two processes yields 
$$\{(g(t)-x_i,l_i) : t_i\leq t, \  x_i<g(t)\} \quad \t{is independent of} \quad
\{(x_i,l_i) : t_i\leq t, \  g(t)\leq  x_i\leq 0\}.$$
Adding that   $\{(x_i,l_i) : t_i\leq t, \  x_i> 0\}$ is
independent of $\{(x_i,l_i) : t_i\leq t, \  x_i\leq  0\}$ and
$g(t)$ is $\{(x_i,l_i) : t_i\leq t, \  x_i\leq  0\}$ measurable,
we get
$$\{(g(t)-x_i,l_i) : t_i\leq t, \  x_i<g(t)\} \quad \t{is independent of} \quad \{(x_i,l_i) : t_i\leq t, \  x_i\geq g(t)\}. \quad$$

We now extend the preceding by incorporating the times $(t_i)_{i\in\N}$. In this direction, we recall that
if $(\w{x}_i,\w{l}_i)_{i\in\N}$  is a  PPP on $\RRR\times \RRR_+$   with intensity $t$d$x\otimes \nu($d$l)$ and 
$(\w{t}_i)_{i\in\N}$ is an iid sequence distributed uniformly on $[0,t]$, then 
$\{(\w{t}_i,\w{x}_i,\w{l}_i) : i \in \N\}$ is distributed as   $\{(t_i,x_i,l_i) : i\in\N, t_i\leq t\}$.
Adding that $g(t)$ is $\{(x_i,l_i) : i\in\N, t_i\leq t\}$ measurable, we get
$$\{(t_i,g(t)-x_i,l_i) : t_i\leq t, \  x_i<g(t)\} \quad \t{is independent of} \quad \{(t_i,x_i,l_i) : t_i\leq t, \  x_i\geq g(t)\}.$$
This ensures that $P_{g(t)}(t)$ is independent of  $(P_{g(t)}^{d(t)}(t),P^{d(t)} (t))$. \\

One can prove similarly that $P^{d(t)} (t)$  is independent of
$(P_{g(t)}(t),P_{g(t)}^{d(t)}(t))$ using that
$(Y^{(t)}_{d(t)+x}-Y^{(t)}_{d(t)})_{x\geq 0}$ is independent of
$(Y^{(t)}_x)_{x\leq d(t)}$ or Lemma 2 in \cite{vb}.
\end{proof}
$\newline$

This guarantees the absence of memory at the left of  $\B (t)$. First we have :
\begin{Cor}
\label{ind}
$(g(t))_{t \in [0,\emph{1/m}]}$   has decreasing càdlàg paths with independent increments.
\end{Cor}
\begin{proof}
Let $0\leq t<t+s\leq \t{1/m}$. The increment $g(t+s)-g(t)$ just
depends on $\Rmt$ and the  point process of files which
arrive after time $t$ at the left of $\B(t)$
$\big\{(t_i,x_i-g(t),l_i) :  t_i>t, x_i<g(t) \big\}$. By the Poissonian property, these two quantities  are
independent and $(g(u) : u \in [0,t])$  is independent  of this point process of files. Moreover
$(g(u) : u \in [0,t])$ is also independent of $(g(t)-\R(t))\cap [0,\infty[$
 by Proposition \ref{indd}. So  $(g(u) : u \in [0,t])$ is independent of
$g(t+s)-g(t)$.
\end{proof}
This explains the observation made in \cite{vb}  Section 3  that the distribution of $g(t)$ is infinitively
divisible (see \cite{Fel} on  page 174 or \cite{Sato} on page 47 for details).

$\newline$
\section{Evolution of the left extremity}

Now we describe the process $(g(t))_{t \in [0,\t{1/m}[}$. We know
that its increments are independent and (\ref{gg})
specifies its marginals. We shall determine its Lévy measure and
prove that its mass is  finite (see \cite{Sato} for terminology). This means 
that the instants when a file arrives at
the left of $\B$ and joins this data block during its storage
do not accumulate before time 1/m, even if $\bar{\nu}(0)=\infty$
(files arrive densely near the data block). Proposition 3 in $\cite{vb}$
ensures that the first time $T_1$ when $0$ is covered, which is also
the first jump time of $(g(t))_{t \in [0,\t{1/m}]}$, is
uniformly distributed on $[0,\t{1/m}]$. Actually the  second jump
time is uniformly distributed in $[T_1,\t{1/m}]$ and  so on ... More
precisely, we have : $\newline$
\begin{Thm}
\label{g}
The jump times of $(g(t))_{t \in [0,\emph{1/m}]}$ are given by an
increasing sequence $(T_i)_{i\in \N}$ which accumulate at $1/\emph{m}$. More
precisely,  using the convention $T_{0}=0$, it holds that for every  $i\geq 1$,  conditionally on $T_{i-1}=t$,   $T_i$  
is independent of $(T_j)_{0\leq j\leq i-1}$   and is uniformly distributed on
 $[t,\emph{1/m}]$.

Then, denoting by $-G_i$ the jump of $(g(t))_{t \in [0,\emph{1/m}]}$ at time $T_i$ for every $i \in \N$, we have
$$g(t):=-\sum_{T_i\leq t} G_i $$
where  $\{(T_i,G_i) : i \in \N \}$ is a PPP  on $[0,\emph{1/m}[\times \RRR^+$ with intensity
$$\emph{d}t\emph{d}x\int_0^{\infty} \P(Y^{(t)}_x \in -\emph{d}l)\bar{\nu} (l).$$

In other words, $(g(t))_{t \in [0,\emph{1/m}]}$ is an additive process and its generating triplet is
$$\left(0,\int_0^{t}\emph{d}s\int_0^{\infty}\P(Y^{(s)}_x \in -\t{d}l) \bar \nu (l),0\right).$$
\end{Thm}
$\newline$
In particular, the interarrival times of $\{T_i : i \in \N\}$ form a 'continuous uniform stick breaking sequence'
(see the residual allocation model in \cite{Pit} on pages 63-64) : the distribution
of $\big((T_{i+1}-T_i)/\t{m}\big)_{i\in\N}$ is the Griffiths-Engen-McCloskey distribution with parameter $(0,1)$
(i.e.   rearranging these increments in the decreasing order yield
 the Poisson-Dirichlet distribution of parameter $(0,1)$). \\ \\
Further,  for every  $i \in \N$, conditionally on  $T_i=t$, the law of $G_i$ is given by
\be
\label{gcondt}
\P(G_i\in \t{d}x)=\t{d}x\frac{1-\t{m}t}{\t{m}}\int_0^{\infty}\P(Y^{(t)}_x \in -\t{d}l) \bar{\nu }(l),
\ee
and as a consequence,
$$\E(G_i)= \big( \frac{1}{(1-\t{m}t)^2}+\frac{1}{2}\frac{\t{m}}{1-\t{m}t} \big)\int_0^{\infty}l^2\nu(\t{d}l). $$
$\newline$
\begin{ex} For the basic example ($\nu=\delta_1$),  conditionally on  $T_i=t$, we have,
$$  \P(G_i\in \t{d}x )=(1-t)e^{-tx}\frac{(tx)^{[x]}}{[x]!}\t{d}x,$$
writing $[x]=\t{sup}\{ n \in \N : n \leq x\}$ and using (\ref{diracY}).
\end{ex}
$\newline$

For the proof, we need the following identity
\begin{Lem}
\label{cv} Let $(S_t)_{t\geq 0}$ be a subordinator with no drift  and  Lévy tail $\bar{\mu}$. Then for all $(t,x) \in \RRR_+^2 $, we have
$$\P(S_t> x)= \int_0^{t} \t{d}s\int_0^x\P(S_s\in \emph{d}b)\bar{\mu}(x-b).$$
\end{Lem}
\begin{proof} As $S$ has no drift,  we have for all $t>0$ and $x>0$,
$$S_t>x\quad \Leftrightarrow \quad \exists ! \ s\in]0,t] \ : \ S_{s^-}\leq x, \ \Delta S_s>x-S_{s^-} \qquad \rm{a.s.}$$
We get then, using also the compensation formula (see \cite{lev} on page 7),
$$\P(S_t>x)=\E(\sum_{0<s\leq t} \ind_{\{S_{s^-}\leq x\}}\ind_{\{\Delta S_s>x-S_{s^-}\}})=\E(\int_0^t\t{d}s\ind_{\{S_s\leq x\}}\bar{\mu}(x-S_s))$$
which completes the proof. One can
also give an analytic proof by computing the Laplace transform of the right hand side for $ q> 0$ and using Fubini :
\bea
&&\int_0^{\infty}\t{d}xe^{-qx}\int_0^{t} \t{d}s\int_0^x \P(S_s\in \t{d}b)\bar{\mu}(x-b)  \nonumber \\
&=& \int_0^{t} \t{d}s  \int_0^{\infty}\mu(\t{d}y) \int_0^{\infty}\P(S_s\in \t{d}b) \frac{e^{-qb}-e^{-q(b+y)}}{q}
= \int_0^{t} \t{d}s  e^{-\phi(q)s}\int_0^{\infty}\mu(\t{d}y)  \frac{1-e^{-qy}}{q} \nonumber \\
&=& \frac{1-e^{-\phi(q)t}}{\phi(q)}\times \frac{\phi(q)}{q}
= \int_0^{\infty}\t{d}xe^{-qx} \P(S_t>x)  \nonumber
\eea
which proves the lemma.
\end{proof}
$\newline$

We are now able to establish Theorem \ref{g}.
\begin{proof}

We know from Corollary \ref{ind} that $(g(t))_{t \in [0, \t{1/m} ]}$ is an additive process. Moreover for every $x\geq 0$,
$(Y^{(t)}_x+x)_{t\geq 0}$ is a subordinator with no drift and Lévy measure $x\nu$ (see ($\ref{defy}$)). So  Lemma \ref{cv} ensures that
\Bea
\P(Y^{(t)}_x>0)&=&\P(Y^{(t)}_x+x>x) \\
&=&\int_0^t\t{d}s \int_0^{x}\P(Y^{(s)}_x+x \in \t{d}b)  x\bar{\nu}(x-b)\\
&=& \int_0^{t}\t{d}s\int_0^{\infty}\P(Y^{(s)}_x \in -\t{d}l) x\bar \nu (l).
\Eea
Using (\ref{gg}), we  get 
$$\E\big(\t{exp}\big(\lambda g(t) \big) \big)=\t{exp}\bigg(\int_0^{\infty} \t{d}x(e^{-\lambda x}-1) 
\int_0^{t}\t{d}s\int_0^{\infty}\P(Y^{(s)}_x \in -\t{d}l) \bar \nu (l)\bigg).$$
So $(g(t))_{t \in [0, \t{1/m} ]}$ is an  additive
process with generating triplet
$$\left(0,\int_0^{t}\t{d}s\int_0^{\infty}\P(Y^{(s)}_x \in -\t{d}l) \bar \nu (l) ,0\right)$$
using   Definition 8.2 and Theorem 9.8 in \cite{Sato}. This characterizes the distribution
of $(g(t))_{t \in [0, \t{1/m} ]}$ (by Theorem 9.8 in \cite{Sato}) and proves that $\{(T_i,G_i) : i \in \N \}$
is a PPP  on $[0,\t{1/m}[\times \RRR^+$ with intensity $\t{d}t\t{d}x\int_0^{\infty} \P(Y^{(t)}_x \in -\t{d}l)\bar{\nu} (l)$.
One can also compute the distribution of $g(t+s)-g(t)$ using the independence of increments and (\ref{gg}) : this proves that
that $g(.)$ is the sum of jumps given by a PPP. \\ \\

By projection, $\{T_i : i\in \N\}$ is a PPP on $[0,1/\t{m}[$ with intensity $\t{m}(1-\t{m}t)^{-1}\t{d}t$. Indeed, for
every $t\in [0,1/\t{m}[$,\\
\bea
\int_0^{\infty}\t{d}x \int_0^{\infty}\P(Y^{(t)}_x \in -\t{d}l) \bar{\nu} (l)&=&\int_0^{\infty}\P(\taup^{(t)}_l \in \t{d}x) \int_0^{\infty} \t{d}l \frac{x}{l} \bar {\nu} (l) \ \ \ \t{using} \ (\ref{egal}) \nonumber \\
&=& \int_0^{\infty} \t{d}l \frac{\E(\taup^{(t)}_l)\bar {\nu} (l)}{l} \nonumber \\
&=& \E(\taup^{(t)}_1)\int_0^{\infty}\bar {\nu} (l)\t{d}l \nonumber \\
&=& \frac{\t{m}}{1-\t{m}t} \ \ \ \ \t{using} \ (\ref{rel}).\nonumber
\eea
Thus, writing $N_t^{t'}:=\t{card}\{ i \in \N : T_i \in ]t,t'] \}$, we have $N_0^{t}<\infty$ a.s. for every $t\in [0,1/\t{m}[$. We
 we can then  sort the times $T_i$ and we have
$$\P(T_{i+1} > t' \mid T_i=t)=\P( N_t^{t'}=0)=\t{exp}\big(-\int_t^{t'} \t{d}s \frac{\t{m}}{1-\t{m}s}\big)=\frac{1-\t{m}t'}{1-\t{m}t},$$
meaning that $T_{i+1}$ is uniformly distributed in $[T_i,\t{1/m}]$. The independence is a consequence
of the Poissonian property of $\{T_i : i\in \N\}$ and we get the theorem. \\

Finally, this proves $(\ref{gcondt})$ and for every $i\in \N$, conditionally on $T_i=t$, we get
\bea
\E(G_i)&=&\frac{1-\t{m}t}{\t{m}}\int_0^{\infty} \t{d}l \frac{\E([\taup^{(t)}_l]^2)\bar \nu (l)}{l} \ \ \ \ \t{using again} \ (\ref{egal}) \nonumber \\
&=&\frac{1-\t{m}t}{\t{m}}\int_0^{\infty} \t{d}l \bar \nu (l)\big(l\big(\frac{\t{m}}{1-\t{m}t}\big)^2+\frac{\int_0^{\infty}l^2\nu(\t{d}l)}{(1-\t{m}t)^3}\big)  \nonumber
\eea
since $[\kappa^{(t)}]'(0)$ is given by (\ref{rel}) and $[\kappa^{(t)}]''(0)$ is given by Proposition 4 in \cite{vb}. 
\end{proof}
$\newline$

\section{The process of remaining data}
We still  consider the files which arrive at the left of $\B$,  the block containing $0$, and   cannot be
entirely  stored at the left of
this block (see Figure \ref{figg}). Such events occur at the  jump times of $(g(t))_{t \in [0,\t{1/m}]}$,  that is at time $T_i$.  We focus here
on the portions  of these files  which cannot be stored at the left of $\B$
and are   shifted  to the right of $\B (T_i-)$  to find a free space. They  are called remaining data
and  denoted by $R_i$. Thus $R_i$
is the quantity of data which arrives at the left of $\B$ at time $T_i$ and is stored at the right of $\B$. Then it is also
 the quantity of data
over $g(T_{i-1}-)$   at  time $T_i$ (see Section 2.1 in \cite{vb} for details) and it is given by
$$ \forall i\geq 1, \ \ \ \ \ R_i:=Y^{(T_i)}_{g(T_{i-1}-)}-I^{(T_i)}_{g(T_{i-1}-)}. $$
$\newline$
We aim at determining the distribution of $\{(T_i,G_i,R_i)\  : \ i \in \N\}$ which is the key to the characterization
of the jumps of $(g(t),d(t))_{t\in [0,\t{1/m}]}$. In that view, we need to describe the arrival of files
which induce the jumps $(G_i,R_i)$. So we consider the  half hardware at the left of $g(t)$, which we turn over, so
that it  is now identified with $\RRR^+$
and its free space is given by $\Rmt$ (see Section 2). The size of free space  and the first free plots
of  this half hardware  are given by the processes $(L^{(t)}_x)_{x\geq 0}$ and  $(D_x^{(t)})_{x\geq 0}$ defined by
$$\forall t \in [0,\t{1/m}[, \ \forall  x\geq 0,  \ \ \ \ \ L^{(t)}_x=\mid \Rmt\cap [0,x]\mid,  \ \ \ \ D_x^{(t)}=\t{inf}\{y>x : y \in \Rmt\}. \ \ \ \ \ \ $$
When at time $t$, a file of length $l$ arrives at location $-x+g(t-)$ on the  hardware (i.e. at location $x$  on the half hardware), it
yields a jump of  $g(.)$   if the free space $L^{(t-)}_x$ between $-x+g(t-)$ and $g(t-)$ is less than $l$. Then the quantity of remaining
data  is  $l-L^{(t-)}_x$  and the jump of the left extremity is $D_x^{(t-)}$ (see Figure \ref{figg}). So we naturally introduce the measure
$\rho^{(t)}$ on $\RRR_+^2$ defined by
$$\rho^{(t)}(\t{d}y\t{d}z):=\int_0^{\infty} \t{d}x \int_0^{\infty} \nu(\t{d}l)  \P(D_x^{(t)} \in \t{d}y,l-L^{(t)}_x \in \t{d}z).$$
In forthcoming Lemma \ref{calc}, we give a useful  alternative expression of $\rho^{(t)}$. This measure gives the intensity of the
point process $\{(T_i,G_i,R_i) : i \in \N\}$, as stated by the following result.


$\newline$
\begin{Thm} $\label{fond}$ $\{(T_i,G_i,R_i) : i \in \N\}$ is a PPP on $[0,\emph{1/m}[\times \RRR_+^2$ with intensity \ $\emph{d}t\rho^{(t)}(\emph{d}y\emph{d}z)$.\\
\end{Thm}

A remarkable consequence is that $(R_i)_{i \in \N}$ is an iid sequence : whereas the rate at which jumps occur
increases as time
gets closer to $\t{1/m}$, the quantity of remaining data keeps the same distribution.
\begin{Cor}
\label{R}
$\{(T_i,R_i) : i \in \N\}$ is a PPP on $[0,\emph{1/m}[\times \RRR^+$ with intensity $\emph{d}t\emph{d}z\frac{\bar{\nu}(z)}{1-\emph{m}t}$. \\
In other words, $(R_i)_{i \in \N}$ is iid,  independent of $(T_i)_{i \in  \N}$ and its distribution is given by :
$$  \P(R_i \in \emph{d}z)=\emph{m}^{-1} \bar{\nu}(z)\emph{d}z, \quad  z\geq 0.$$
\end{Cor}
$\newline$
\begin{ex}
Using the expression of $\rho^{(t)}$ given by Lemma  $\ref{calc}$ below, the 
expressions $(23)$ and $(24)$ in \cite{vb} yield an expression of $\rho^{(t)}$ for
the basic example and the gamma distribution which is quite heavy and not mentioned here. Nonetheless
the quantity of remaining data can be often calculated
explicitly. For the basic example ($\nu=\delta_1$), the
remaining data are uniform random variables on $[0,1]$. For the
exponential distribution ($\nu(\t{d}l)=\ind _{\{l\geq 0\}}e^{-l}\t{d}l$),
the remaining data  are also exponentially distributed.
\end{ex}
$\newline$

The proofs of these results are organized as follows. \\
First, in Lemma \ref{calc}, we give a more explicit  expression of $\rho^{(t)}$ which will be  useful for the proofs and will
 enable  us to derive
 Corollary \ref{R} from Theorem $\ref{fond}$. \\
 Second, we prove
that $\rho^{(t)}$  gives the intensity of the point process $\{(T_i,G_i,R_i) : i \in \N\}$ (Lemma \ref{lim}). That is
for every $t \in [0,\t{1/m}[$ and $A=]a_1,b_1]\times ]a_2,b_2] \subset \RRR_+^2$, we have :
$$\lim_{h\rightarrow 0} \frac{\P(\exists i \in \N : T_i \in ]t,t+h], \ (G_i,R_i) \in  A)}{h}=\rho^{(t)}(A).$$
The lowerbound appears  naturally by considering the  arrival of one
single file independently  of the past which induces a jump of the
left extremity, as described at the beginning of this section
(see also  Figure \ref{figg}). However, in the case
$\bar{\nu}(0)=\infty$, some jumps of the left extremity could be
due to the successive arrival of many files  during a short time
interval $]t,t+h]$. Thanks to
Theorem \ref{g}, we already know the rate at which jumps occur (i.e. the total intensity). This will give us the upperbound.\\
Finally, we prove that  the point process $\{(T_i,G_i,R_i) : i \in \N\}$ enjoys a
memoryless property (Lemma \ref{indp}), which is a direct consequence of results of Section 3.  We  get then
the complete description of this point process, which enables us to prove Theorem  $\ref{fond}$. Corollary \ref{R} follows
by integrating $\rho^{(t)}$ with respect to the first coordinate. \\
$\newline$

Recall the notation in Theorem \ref{loi} and (\ref{kappat}).
\begin{Lem}\label{calc} For every $t \in [0,\emph{1/m}[$, the measure $\rho^{(t)}(\emph{d}y\emph{d}z)$ can also be expressed as 
\setlength\arraycolsep{0pt}
\bea
& &\emph{d}z\int_z^{\infty} \nu(\emph{d}l) \bigg(\P(\taum^{(t)}_{l-z} \in \emph{d}y)+\int_0^{y}\P(\taum^{(t)}_{l-z} \in \emph{d}x)(y-x)\Pi^{(t)}(\emph{d}y-x)\bigg) \nonumber \\
&=& \int_z^{\infty} \nu(\emph{d}l)(l-z) \bigg(y^{-1} \emph{d}y\P(Y^{(t)}_y+l \in \emph{d}z)+\int_0^{y}\P(Y^{(t)}_x+l \in  \emph{d}z)(yx^{-1}-1)\Pi^{(t)}(\emph{d}y-x)\bigg) \nonumber
\eea
\end{Lem}
\begin{proof} By Lemma 1.11 in Chapter 1 of $\cite{mono}$ applied to  $(\taum^{(t)}_{x})_{x\geq 0}$, 
 we have for all $a,b\geq 0$ and $q>0$ ($t$ is fixed and omitted in the notation),
$$\int_0^{\infty}\t{d}xe^{-qx}\E(\t{exp}(-bL_x-aD_x))=\frac{\kappa(a+q)-\kappa(a)}{q(\kappa(a+q)+b)}.$$
Letting $q\rightarrow 0$, we get
$$\int_0^{\infty} \t{d}x \E(\t{exp}(-bL_x-aD_x))=\frac{\kappa'(a)}{\kappa(a)+b}=\int_0^{\infty}\t{d}ze^{-bz}\kappa'(a)e^{-\kappa(a)z}.$$
From $\kappa'(a)=\int_0^{\infty}e^{-ay}(\delta_0(\t{d}y)+y\Pi(\t{d}y))$ and $e^{-\kappa(a)z}=\int_0^{\infty}e^{-ay}\P(\taum_z\in \t{d}y)$, we deduce
\be
\label{convo}
\int_0^{\infty} \t{d}x \E(\t{exp}(-bL_x-aD_x))=\int_0^{\infty}\t{d}z\int_0^{\infty}\gamma_z(\t{d}y)e^{-bz-ay},
\ee
where $\gamma_z$ is the convolution of $\delta_0(\t{d}y)+y\Pi(\t{d}y)$ and $\P(\taum_z\in \t{d}y)$. Thus,
\Bea
\gamma_z(\t{d}y)&=& \int_0^{y}\P(\taum_z\in \t{d}x)(\delta_0(\t{d}y-x)+(y-x)\Pi(\t{d}y-x))\nonumber \\
&=&\P(\taum_z\in \t{d}y)+\int_0^{y}\P(\taum_z\in \t{d}x)(y-x)\Pi(\t{d}y-x). \nonumber
\Eea
And the identification of Laplace transforms in ($\ref{convo}$) entails that
\be
\label{calclaplace}
\int_0^{\infty} \t{d}x \P(L_x \in \t{d}z, D_x \in \t{d}y)=\t{d}z\big(\P(\taum_z\in \t{d}y)+\int_0^{y}\P(\taum_z\in \t{d}x)(y-x)\Pi(\t{d}y-x)\big),
\ee
which proves the first identity  of the lemma integrating  with respect to $l$. Using (\ref{egal}) gives  the second one.
\end{proof}
\begin{Rque}
A recent work of Winkel (Theorem 1 in $\cite{Wink}$) enables to calculate differently the law of
$\P(L_x \in \t{d}z, D_x \in \t{d}y)$ ($L_x$ corresponds  to $T_x$ in $\cite{Wink}$ and $D_x$  to $X(T_{x-})+\Delta_x$) :
$$\int_0^{\infty} \t{d}x \P(L_x \in \t{d}z, D_x \in \t{d}y)=\t{d}y\P(H_y\in \t{d}z)+\t{d}z\int_0^{\infty}\P(\taum_x\in \t{d}x)(y-x)\Pi(\t{d}y-x),$$
where $H_x=\t{inf}\{a\geq 0, \taum_a=x \}$. Then observe that the  measures on $\RRR_+^2$  $\t{d}y\P(H_y\in \t{d}z)$ and
$\t{d}z\P(\taum_z\in \t{d}y)$ coincide by computing their Laplace transform using (4) in
$\cite{Wink}$. This proves  $(\ref{calclaplace})$.
\end{Rque}
$\newline$

Second, for every  Borel set $B$ of
$[0,\t{1/m}[\times \RRR^{2}_+$, we define   $N_{B}:=\t{card}\{i \in \N : (T_i,G_i,R_i) \in B \}$ and we say that $A$ is
 a  rectangle of $D\subset \RRR^{\t{d}}$ if $A$ is a subset of $D$  of the form
$$\{x=(x_1,x_2,..,x_{\t{d}}), a_1<x_1\leq b_1,..,a_{\t{d}}<x_{\t{d}}\leq b_{\t{d}}\}.$$
Then, we have
\begin{Lem}
\label{lim}
For all $t \in [0,\emph{1/m}[$ and $A$ rectangle  of $\RRR_+^2$, we have :
$$\lim_{h\rightarrow 0} \frac{\P(N_{]t,t+h]\times A}\geq 1)}{h}=\rho^{(t)}(A).$$
\end{Lem}
\begin{proof} First we prove the lowerbound. Second, we check that  the convergence
holds for $A=\RRR_+^2$. \\ \\
$\bullet$
Let $\epsilon >0$, $A=]a,b]\times ]c,d]$  and work conditionally on  $\Rmt$. We consider   a file
labelled $i$ which arrives at time
$t_i \in ]t,t+h]$ at location $x_i<g(t)$. We put $\w{x_i}:=g(t)-x_i \geq 0$
the arrival point on the half line at the left of $g(t)$ and require that
$$l_i-L^{(t)}_{\w{x_i}} \in ]c,d-\epsilon],  \ D_{\w{x_i}}^{(t)} \in ]a,b-\epsilon],  \  \vert L^{(t_i-)}_{\w{x_i}} -L^{(t)}_{\w{x_i}}\vert\leq
\epsilon , \  \vert D_{b}^{(t_i-)} -D_{b}^{(t)} \vert \leq \epsilon. $$
Then file $i$ verifies
$$\ l_i-L^{(t_i-)}_{\w{x_i}} \in ]c,d],  \ D_{\w{x_i}}^{(t_i-)} \in ]a,b].$$
So this file induces a jump of the left extremity  and  $N_{]t,t+h]\times A}\geq 1$ (see the beginning of this section or Figure \ref{figg} for details) and we get the lowerbound :
\setlength\arraycolsep{0pt}
\bea
&&\P\big(N_{]t,t+h]\times A}\geq 1 \ \mid  \ \Rmt \big)     \nonumber \\
&& \ \ \ \ \ \ \ \geq  \ \P\big(\exists i \in \N\  : \   t_i \in ]t,t+h], \ l_i-L^{(t)}_{\w{x_i}} \in ]c,d-\epsilon], \ D_{\w{x_i}}^{(t)} \in ]a,b-\epsilon], \nonumber \\
&& \ \ \ \ \ \ \ \ \ \ \ \  \  \vert L^{(t_i-)}_{\w{x_i}} -L^{(t)}_{\w{x_i}}\vert \leq \epsilon , \  \vert D_{b}^{(t_i-)} -D_{b}^{(t)} \vert \leq \epsilon \ \mid  \ \Rmt \big) \nonumber \\
\label{cvint}
&& \ \ \ \  \ \ \ \geq \ A_t(h) . B_t(h)
\eea
where
\bea
&& A_t(h):=\P\big(\exists i \in \N :  \  t_i \in ]t,t+h], l_i-L^{(t)}_{\w{x_i}} \in ]c,d-\epsilon], D_{\w{x_i}}^{(t)} \in ]a,b-\epsilon]  \ \mid  \ \Rmt \big), \nonumber \\
&& B_t(h):=\P\big( \sup_{t'\in[t,t+h]}\{ \vert L_b^{(t')}-L_b^{(t)}\vert \} \leq \epsilon, \sup_{t'\in[t,t+h]}\{\vert D_{b}^{(t')}-D_{b}^{(t)} \vert\} \leq \epsilon  \ \mid  \ \Rmt \big). \nonumber
\eea
1) By Theorem $\ref{g}$, $\P(N_t^{t+h} \ne 0) \stackrel{h\rightarrow 0}{\longrightarrow }0$ so a.s for $h$ small
enough, $g(t+h)=g(t)$. Then,  using the Hausdorff metric on $\RRR_+$  (denoted by $\H(\RRR_+)$ in Section 2 in  $\cite{vb}$), we have
$$\overset{\longleftarrow }{\scriptstyle{ \mathcal{R}(t+h)}} \stackrel{h\rightarrow 0}{\longrightarrow } \Rmt \quad \t{a.s.}$$
Then $B_t(h)$ converges a.s. to $1$ as $h$ tends to $0$. \\
$\newline$
2) As $\{(t_i,\w{x_i},l_i) : i\in\N ,t_i \in ]t,t+h],  x_i<g(t)\}$ is a PPP on $]t,t+h]\times \RRR_+^2$ with intensity d$t\otimes $d$x\otimes \nu$(\t{d}$l$)  independent of $\Rmt$,

$$A_t(h)=1-\t{exp}\big(-h\int_0^{\infty} \t{d}x \int_0^{\infty} \nu(\t{d}l)\ind_{\{l-L^{(t)}_x \in ]c,d-\epsilon],D_x^{(t)} \in ]a,b-\epsilon]\}}\big) \qquad \t{a.s.}$$
This term is a.s. equivalent when $h$ tends to $0$ to
$$h\int_0^{\infty} \t{d}x \int_0^{\infty} \nu(\t{d}l)\ind_{\{l-L^{(t)}_x \in ]c,d-\epsilon],D_x^{(t)} \in ]a,b-\epsilon]\}}.$$
$\newline$

Then, letting $h\rightarrow 0$ in $(\ref{cvint})$, 1) and 2) give
$$\liminf_{h\rightarrow 0} \frac{\P\big(N_{]t,t+h]\times A}\geq 1 \ \mid  \ \Rmt \big)}{h}\geq \int_0^{\infty} \t{d}x \int_0^{\infty} \nu(\t{d}l)\ind_{\{l-L^{(t)}_x \in ]c,d-\epsilon],D_x^{(t)} \in ]a,b-\epsilon]\}} \ \ \t{a.s.}$$
Integrating this inequality  and using Fatou's lemma yield
\bea
\liminf_{h\rightarrow 0}\frac{\P\big(N_{]t,t+h]\times A}\geq 1\big)}{h}&\geq&  \E\big(\int_0^{\infty} \t{d}x\int_0^{\infty}\nu(\t{d}l)\ind_{\{l-L^{(t)}_x \in ]c,d-\epsilon],D_x^{(t)} \in ]a,b-\epsilon]\}}\big) \nonumber \\
&\geq & \rho^{(t)}(]a,b-\epsilon]\times ]c,d-\epsilon]). \nonumber
\eea
As $\rho^{(t)}(]a,b]\times \{d\}\cup \{b\}\times ]c,d])=0$  (use the two equalities of Lemma \ref{calc}), we get letting $\epsilon$ tend to $0$  :  
$$\liminf_{h\rightarrow 0}\frac{\P\big(N_{]t,t+h]\times A}\geq 1\big)}{h}\geq \rho^{(t)}(A).$$ \\ \\
$\bullet$ 
We derive the upperbound from Theorem \ref{g}. First,
$$\frac{\P\big(N_{]t,t+h]\times \RRR_+^2}\geq 1\big)}{h}=\frac{\P(\exists i \in \N   \ : \ T_i \in ]t,t+h])}{h}\stackrel{h\rightarrow 0}{\longrightarrow } \frac{\t{m}}{1-\t{m}t}.$$
and identity $(\ref{calcint})$ below gives
$$\rho^{(t)}(\RRR_+^2)=\frac{\t{m}}{1-\t{m}t}.$$
So we just need to prove the following result : Let $(\mu_n)_{n\in \N}$ and $\mu$ be finite measures on $\RRR_+^2$ such that  for every  $A$ rectangle of $\RRR_+^2$ :
$\liminf_{n \rightarrow \infty} \mu_n(A)\geq \mu(A)$ \ and \ $\lim_{n\rightarrow \infty} \mu_n(\RRR_+^2)=\mu(\RRR_+^2)$. Then
for every $A$ rectangle of $\RRR_+^2$, $\lim_{n\rightarrow \infty} \mu_n(A)= \mu(A)$. \\
In that view, suppose there exist a rectangle $A$, $\epsilon>0$  and a sequence of integers $k_n$  such that $\mu_{k_n}(A)\geq \mu(A)+\epsilon$. Choose $B$ union of disjoint rectangles all disjoint from $A$ such that
 $\mu(B\cup A)\geq \mu(\RRR_+^2)-\epsilon/2$. Then,
$$\liminf_{n\rightarrow \infty} \mu_{k_n}(\RRR_+^2)\geq \liminf_{n\rightarrow \infty} \mu_{k_n}(A\cup B)\geq \mu(A)+\epsilon+\mu(B)\geq \mu(\RRR_+^2)+\epsilon /2,$$
which is a contradiction with $\lim_{n\rightarrow \infty} \mu_n(\RRR_+^2)=\mu(\RRR_+^2)$.
\end{proof}
$\newline$

To prove the theorem, it  remains to prove the absence of memory.
\begin{Lem}
\label{indp}
Let $t\in[0,\emph{1/m}[$, then $\big\{(T_i,G_i,R_i)\  : \ i \in \N,\ T_i\leq t\big\}$ is
independent of $\big\{(T_i,G_i,R_i)\  : \ i \in \N, \ T_i> t\big\}.$
\end{Lem}
\begin{proof} First $\big\{(T_i,G_i,R_i)\  : \ T_i\leq t\big\}$ is given by  $\big\{(t_i,l_i,x_i) : t_i\leq t, x_i \in [g(t),d(t)]\big\}$.
 Moreover $\big\{(T_i,G_i,R_i)\  : \ T_i> t\big\}$  depends on $(\R(t)-g(t))\cap ]-\infty,0]$  and
$\big\{(t_i,x_i-g(t),l_i)\  : \ t_i>t,\ x_i< g(t)\big\}$ which are independent.
Moreover  $(\R(t)-g(t))\cap ]-\infty,0]$  is independent of $\big\{(t_i,l_i,x_i) : t_i\leq t, x_i \in [g(t),d(t)]\big\}$
by  Lemma $\ref{indpp}$ and so is  $\big\{(t_i,x_i-g(t),l_i)\  : \ t_i>t,\ x_i< g(t)\big\}$ by Poissonian property. This proves the result.
\end{proof}
$\newline$
We can now prove the theorem and its corollary.
\begin{proof}[Proof of  Theorem \ref{fond}]
We prove now  that  for every $B$ finite union of disjoint rectangles of $[0,\t{1/m}[\times \RRR_+^2$:
\be
\label{egmes}
\P(N_B=0)=e^{-\gamma(B)}, \ \ \ \t{where} \ \ \gamma(\t{d}t\t{d}y\t{d}z)=\t{d}t\rho^{(t)}(\t{d}y\t{d}z).
\ee
As $\gamma$ is non atomic (use Lemma $\ref{calc}$), this will ensure that $\{(T_i,G_i,R_i) : i \in \N\}$ is a PPP
 with intensity  $\gamma$  (use Renyi's Theorem $\cite{King}$). \\

Let $t \in [0,\t{1/m}[$ and $A$ a finite union of rectangles of $\RRR_+^2$. We consider $H(s):=\P(N_{]t,t+s]\times A}=0)$ for $s \in [0,\t{1/m}-t[$. Lemma $\ref{indp}$
entails that $$H(s+h)=\P(N_{]t,t+s]\times A}=0)\P(N_{]t+s,t+s+h]\times A}=0)=H(s)\P(N_{]t+s,t+s+h]\times A}=0).$$
We write $A=\sqcup_{i=1}^{N} A_i$ where $A_i$ rectangle of $\RRR_+^2$. Theorem \ref{g} and Lemma \ref{lim} ensure respectively that for all $1\leq i,j\leq N$ such that $i \ne j$:
$$\lim_{h\rightarrow 0} \frac{\P(N_{]t,t+h]\times A_i}\geq 1, \ N_{]t,t+h]\times A_j}\geq 1)}{h}=0 \ \ \ ; \ \ \ \lim_{h\rightarrow 0} \frac{\P(N_{]t,t+h]\times A_i}\geq 1)}{h}=\rho^{(t)}(A_i).$$
Then
$$\lim_{h\rightarrow 0} \frac{\P(N_{]t,t+h]\times A}\geq 1)}{h}=\sum_{i=1}^{N} \lim_{h\rightarrow 0} \frac{\P(N_{]t,t+h]\times A_i}\geq 1)}{h}=\rho^{(t)}(A),$$
and the derivative of $H$ is given by
$$\lim_{h\rightarrow 0} \frac{H(s+h)-H(s)}{h}= H(s)\lim_{h\rightarrow 0} \frac{1-\P(N_{]t+s,t+s+h]\times A}=0)}{h}= H(s)\rho^{(t+s)}(A).$$
Thus $H(s)$ satisfies a differential equation of order $1$ and we get (\ref{egmes}) for $B=]t,t+s]\times A$.
$$H(s)=\t{exp}\big(-\int_0^s \t{d}u\rho^{(t+u)}(A) \big)=\t{exp}\big(-\int_t^{t+s} \t{d}u\rho^{(u)}(A) \big)=e^{-\gamma(]t,t+s]\times A)}$$
Using again Lemma $\ref{indp}$ and additivity of measures proves (\ref{egmes}) for every $B$ finite union of rectangles of $[0,\t{1/m}[\times \RRR^+\times\RRR^+$.
\end{proof}
$\newline$
\begin{proof}[Proof of  Corollary \ref{R}] As  projection of the  PPP $\{(T_i,G_i,R_i) : i \in \N\}$, $\{(T_i,R_i) : i \in N \}$ is a PPP with
intensity d$t\int_{y\in[0,\infty]}\rho^{(t)}(\t{d}y\t{d}z)$. By Lemma $\ref{calc}$, we have :
\bea
\int_{y\in[0,\infty]}\rho^{(t)}(\t{d}y\t{d}z) &=& \t{d}z\big(\bar{\nu}(z)+\int_z^{\infty}\nu(\t{d}l)\int_0^{\infty}\P(\taum^{(t)}_{l-z} \in \t{d}x)\int_x^{\infty} \Pi^{(t)}(\emph{d}y-x)(y-x)\big)   \nonumber \\
&=& \t{d}z\bar{\nu}(z)(1+\int_0^{\infty}\Pi(\t{d}y)y) \nonumber \\
\label{calcint}
&=& \t{d}z\frac{\bar{\nu}(z)}{1-\t{m}t}  \ \ \ \ \ \ \ \ \t{by}  \ (\ref{rel})
\eea
which gives the intensity of $\{(T_i,R_i) : i \in \N \}$. In other words,  $(R_i)_{i\in\N}$
is an iid sequence independent of $(T_i)_{i\in\N}$ such that
$\P(R_i \in \emph{d}z)=\t{m}^{-1} \bar{\nu}(z)\t{d}z, \ (z\geq 0)$.
\end{proof}

$\newline$

\section{Evolution of  $\B$}

$\newline$ 
The processes $(g(t))_{t \in [0,\t{1/m}]}$ and $(d(t))_{t \in
[0,\t{1/m}]}$  of the left and the right extremities of $\B$  have a quite different evolution, even though their
one-dimensional distributions coincide. The process
$(d(t))_{t \in [0,\t{1/m}]}$ jumps each time $(g(t))_{t \in
[0,\t{1/m}]}$ jumps  and each time  a file arrives on $\B$. More
precisely, there are two kinds of jumps
of $(\B(t))_{t \in [0,\t{1/m}]}$ corresponding respectively to : \\
- files which arrive at the left of $\B$ and cannot be entirely stored at its left (recall the previous section). These files
induce the jumps $(-G_i,D_i)$ of
the  extremities of $\B$ at time $T_i$ independently of the past (see Figure \ref{figg}). \\
- files which arrive on $\B$. These files induce  jumps of the right extremity $d(.)$ only,  with  total  rate  equal to
$l(t)\bar{\nu}(0)$ (see Figure \ref{figgg}). This rate is infinite when $\bar{\nu}(0)=\infty$. Observe also that the jumps depend from the past of $\B$ through the value
of the length $l(t)$.\\
Note that  a file which arrives at the left of $\B(t-)$ at time $t$ with remaining data of size $R$
 induces the same jump of the right extremity as a file of size $R$ which arrives on $\B(t-)$ at time $t$.
Obviously, the other files (files which are entirely stored at the left of $\B$ or  which arrive at the right of $\B$) do
not yield a jump of $\B$.\\ \\
Thus,  we define
$$D_i:=d(T_i)-d(T_i^-)$$
 and we decompose the process $(g(t),d(t))_{t\in[0,\t{1/m}[}$ into two processes $(C^1(t))_{t\in[0,\t{1/m}[}$ and
$(C^2(t))_{t\in[0,\t{1/m}[}$, which  give  the variation of the extremities  of $\B$ respectively  at times $(T_i)_{i\in\N}$
(due to the arrival of a file at the
left of $g(t)$) and  between successive times $(T_i)_{i\in\N}$ (due to the arrival of files on $\B(t)$). That is, for every $t\in[0,1/\t{m}[$,
$$C^1(t):=\sum_{T_i\leq t} (-G_i,D_i), \qquad   C^2(t):=\big(0, \sum_{\substack{0\leq s\leq t \\ s \notin \{T_i : i\in \N\}}}
\Delta d (s)\big),$$
$$(g(t),d(t))=C^1(t)+C^2(t).$$
$\newline$

First, we specify  the distribution of $(C^1(t))_{t\in[0,1/\t{m}]}$ (see below for the proofs).
\begin{Pte}
\label{pppd}
The point process $\big\{(T_i,G_i,D_i) : i \in \N \big\}$ is a PPP on $[0,\emph{1/m}]\times \RRR_+^2$ with intensity
$\emph{d}t\mu^{(t)}(\emph{d}y\emph{d}x)$, where
$$\mu^{(t)}(\emph{d}y\emph{d}x)=\int_0^{\infty}\rho^{(t)}(\emph{d}y\emph{d}z)\P(\taup^{(t)}_z\in \emph{d}x).$$
\end{Pte}
$\newline$

We can now specify the distribution of the process $(g(t),d(t))_{t\in[0,\t{1/m}[}$ as follows.
\newpage
\begin{Thm} \label{couple}  $(g(t),d(t))_{t\in[0,\emph{1/m}[}$ is a pure jump Markov process equal to
$(C^1(t)+C^2(t))_{t\in [0,\emph{1/m}]}$ such that for all $0\leq t\leq t+s\leq \emph{1/m}$,  \\ \\
(i) $C^1(t+s)-C^1(t)$ is independent of $(g(u),d(u))_{u\in [0,t]}$ .   \\ \\
(ii) Conditionally on  $l(t)=l$, $C^2(t+s)-C^2(t)$ is independent of $(g(u),d(u))_{u\in [0,t]}.$   \\
Conditionally   also on $T_i\leq t\leq t+s<T_{i+1}$ for some $i \in \N$ :
$$C^2(t+s)-C^2(t)\stackrel{d}{=}(0,\taup^{(t+s)}_{S_{sl}}),$$
where $(S_x)_{x\geq 0}$ is a subordinator with no drift and Lévy measure $\nu$, which is independent 
of $(\taup^{(t+s)}_x)_{x\geq 0}$.
\end{Thm}
$\newline$

We recall that vague convergence of measures on $A$ is the convergence of  the integrals of measures against
continuous functions with compact support in $A$. The  jump rate of $(g(t),d(t))_{t \in [0,\t{1/m}[}$ is then given by   :
\begin{Cor} \label{rate}
If $t \in [0,\emph{1/m}[$, we have the following vague convergence of measures on $[0,\infty[\times ]0,\infty[$ when  $h$ tends  to $0$ :
$$h^{-1}\P(g(t)-g(t+h) \in \emph{d}y, \ d(t+h)-d(t) \in \emph{d}x \ \mid \ l(t)=l) \ \stackrel{w}{\Longrightarrow} $$
$$\mu^{(t)}(\emph{d}y\emph{d}x)+l\delta_0(\emph{d}y)\int_0^{\infty}\nu(\emph{d}z)\P(\taup^{(t)}_{z}\in \emph{d}x).$$
\end{Cor}
$\newline$
$\newline$

We begin with two lemmas which state the independences needed for the proofs.

\begin{Lem}
\label{indD}
$\{(T_i,G_i,D_i) : i \in \N, \ T_i>t\}$ is independent of $(g(u),d(u))_{u\in [0,t]}$.
\end{Lem}
\begin{proof}
Using $(\ref{Di})$ below, we see that $\{(T_i,G_i,D_i) : i \in \N, \ T_i>t\}$ is given by
$$\{(T_i,G_i,R_i) : i \in \N, T_i> t\} \ \ \  \t{and}  \ \ \ \ (\Rp(s))_{s>t}. $$
These quantities depend from the past through $(\Rmt,\Rpt)$ which is independent of $(g(u),d(u))_{u\in [0,t]}$
by  Proposition $\ref{indd}$.
\end{proof}
$\newline$
\begin{Lem}
\label{inddd}
Let $i\in \N$ and $0\leq t'<t\leq \emph{1/m}$. Conditionally on  $T_{i-1}=t'$ and $T_i=t$, $(\overset{\longrightarrow }{\scriptstyle{ \mathcal{R}(u)}})_{u\in [t',t[}$   is independent of the point process $P_{g(t')}(t)$.
\end{Lem}
\begin{proof}
Conditioning by $T_{i-1}=t'$ and $T_i=t$ ensures that all the data arrived at the left of $g(t')$ during the time interval $[t',t[$ are 
stored at the left of $g(t')$. So $(\overset{\longrightarrow }{\scriptstyle{ \mathcal{R}(u)}})_{u\in [t',t[}$ depdns only on
the point process $P_{g(t')}^{d(t')}(t)\cup P^{d(t')}(t)$
which is independent of $P_{g(t')}(t)$ by Lemma $\ref{indpp}$.
\end{proof}
$\newline$
\begin{proof}[Proof of  Proposition \ref{pppd}]
At time $T_i$, the quantity of remaining data $R_i$ is  stored  at the right of $\B (T_i-)$. It induces a jump $D_i=d(T_i)-d(T_i-)$
of the right extremity
which is equal  to $R_i$ plus the sum   of the lengths of blocks  at the right of $\B (T_i-)$
which are reached during the storage of these data (see Figure 2). More precisely :
\bea
D_i&=&
\t{inf}\{x\geq 0, \ \mid \R(T_i-) \cap [d(t),d(t)+x[\vert=R_i\} \nonumber \\
&=&\t{inf}\{x\geq 0, \ \mid \overset{\longrightarrow }{\scriptstyle{ \mathcal{R}(T_i-)}}\cap [0,x] \mid = R_i\} \nonumber \\
\label{Di}
&=&\taup^{(T_i-)}_{R_i},
\eea
by definition of $\taup$ (see Section 2). Lemma $\ref{inddd}$ ensures that conditionally on $T_i=t$, $(\taup^{(T_i-)}_x)_{x\geq 0}$ is independent of
$(G_i,R_i)$ and distributed as $(\taup^{(t)}_x)_{x\geq 0}$. Then denoting by  $\mu_t$ the law of $(G_i,D_i)$ conditioned by $T_i=t$, we have
\be
\label{mutild}
\mu_t(\t{d}y\t{d}x)=\P(G_t\in \t{d}y, \ \taup^{(t)}_{R_t} \in \t{d}x),
\ee 
where $(G_t,R_t)$ is a random variable independent of $(\taup^{(t)}_x)_{x\geq 0}$ and distributed as $(G_i,R_i)$ conditioned on $T_i=t$.\\

By Lemma $\ref{indD}$, $\big\{(T_i,G_i,D_i) :  i \in \N, \ T_i> t
\big\}$ is independent of  $\big\{(T_i,G_i,D_i) : i \in \N,
T_i\leq t \big\}$. Then conditionally on   $(T_i)_{ i \in \N}$,
$(G_i,D_i)_{ i \in \N}$ are independent. Adding that $\{T_i : i \in \N\}$ is a PPP on $[0,\t{1/m}]$
with intensity $\t{d}t\t{m}/(1-\t{m}t)$  ensures that
$\big\{(T_i,G_i,D_i) : i \in \N \big\}$ is a (marked) PPP with
intensity 
$$\frac{\t{m}}{1-\t{m}t}\t{d}t\mu_t(\t{d}y\t{d}x).$$
Furher, by  ($\ref{mutild}$), this intensity is eqaul to 
$$ \t{d}t \int_0^{\infty}\P(\taup^{(t)}_z\in \t{d}x)\frac{\t{m} }{1-\t{m}t} \P(G_t \in \t{d}y, \ R_t \in \t{d}z)
= \t{d}t \int_0^{\infty}\P(\taup^{(t)}_z\in \t{d}x)\rho^{(t)}(\t{d}y \t{d}z)
$$
using Theorem \ref{fond}. This completes the proof.
\end{proof}

$\newline$

\begin{proof}[Proof of  Theorem \ref{couple}]
$\newline$
$(i)$ \ Thanks to Lemma $\ref{indD}$, $C^1(t+s)-C^1(t)$ is independent of $(g(u),d(u))_{u\in [0,t]}$.\\  \\
$(ii)$ \ We condition by  $T_i\leq t\leq t+s<T_{i+1}$ for some $i \in \N$ and   $l(t)=l$. Then $g(t+s)-g(t)=0$ and
 no data arrived at the left of $\B (t)$ during the time interval $]t,t+s]$ is stored at the right of this block. So the increment
$d(t+s)-d(t)$ is caused by files arriving on $\B (t)$ : they are stored at the right on $\B (t)$ and may join
data already stored. Note that we can change the order of arrival of files between $t$ and $t+s$
(use identity $(\ref{comp})$). Thus, we first  store  the files which arrive at the right of $d(t)$ between times $t$ and
$t+s$, then the files which arrive on $\B(t)$ between times $t$ and
$t+s$ and we forget the  files which arrive at the left of $g(t)$. \\

STEP 1 : At time $t$, we consider the half hardware  at the right  of $d(t)$ which we identify with $[0,\infty[$. Its  free
space is equal to $\Rpt$. We store the files $i \in \{i \in \N : t_i \in ]t,t+s], x_i >d(t)\}$ on this half hardware
 $[0,\infty[$  at location $x_i-d(t)$ following the process described in Introduction
(the size of the file $i$ is still $l_i$). Following Section 2.1 in \cite{vb}, we get the counterpart of
the characterization of the free space $(\ref{comp})$. That is, the new  free  space
of the half hardware is  equal to  $\{ x \geq 0  :  \w{Y}_{x}=\w{I}_{x}\}$ , where for every $x\geq 0$,
$$ \w{Y}_x=-x+\sum_{\substack{0\leq t_i\leq t+s \\ d(t)\leq x_i\leq d(t)+x}} l_i, \ \
\ \ \w{I}_x:=\t{inf}\{\w{Y}_y : 0\leq y\leq x\}.$$
Using Lemma \ref{indpp}, we see that
$\{(t_i,x_i -d(t),l_i) :  x_i \geq d(t)\}$ is  a   PPP on $\RRR^{+3}$ with intensity d$t\otimes \t{d}x\otimes \nu(\t{d}l)$. Then,
$$\big(\w{Y}_x\big)_{x\geq0}\stackrel{d}{=}\big(Y^{(t+s)}_x\big)_{x\geq0}$$
is a Lévy process with Laplace exponent $\Psi^{(t+s)}$.
As  $[\Psi^{(t+s)}]'(0)<0$,  $\big(\w{Y}_x\big)_{x\geq0}$ is regular
for $]-\infty,0[$, in the sense that it takes negative values for some arbitrarily small $x$ (Proposition 8 on page 84 in  \cite{lev}). So for every stopping time $T$ such that
$\w{Y}_T=\w{I}_T$, there is the identity $T=\inf\{ z\geq 0 : \w{Y}_{z}<\w{Y}_T\}$.  This ensures that
the  free space $\{ x \geq 0  :  \w{Y}_{x}=\w{I}_{x}\}$ of the half hardware  is the range of
$(\w{\tau}_{x})_{x\geq 0}$ defined by
$$\w{\tau}_{x}:=\t{inf}\{z\geq 0 : \w{Y}_z<-x\}.$$
By Theorem 1 on page 189 in  \cite{lev}, $(\w{\tau}_{x})_{x\geq0}$ is a subordinator with Laplace
exponent  $\kappa^{(t+s)}$, which is the inverse function of $-\Psi^{(t+s)}$.
So  $(\w{\tau}_{x})_{x\geq0}$ is distributed as  $(\taup^{(t+s)}_x)_{x\geq 0}$. By Lemma 1 again,  $\{(t_i,x_i -d(t),l_i) :  x_i >d(t)\}$
is independent of $(g(u),d(u))_{u \in [0,t]}$. So  $(\w{\tau}_{x})_{x\geq0}$ is independent of $(g(u),d(u))_{u \in [0,t]}$.
$\newline$

STEP 2 : To obtain the covering $\C(t+s)$, we now store  the files
$\{i : t_i \in ]t,t+s], x_i \in [g(t),d(t)[\}$. It amounts to
store these files  in the first free spaces (i.e.
as much on the left as possible) of the half hardware  considered above, whose
free space is the range of $(\w{\tau}_{x})_{x\geq0}$. The variation of the right
extremity is  equal to the sum of the sizes of these files, say $S_t^{t+s}$, plus
the sizes of the lengths of the blocks of the half hardware joined
during their storage. That is, as for   (\ref{Di}),
$$ C^2(t+s)-C^2(t)=(0,\w{\tau}_{S_t^{t+s}}), \qquad  \t{where} \quad  S_t^{t+s}:= \sum_{\substack{t<t_i\leq  t+s \\ x_i \in [g(t),d(t)[}} l_i.$$
Conditionally on $l(t)=l$, by Poissonian property, $S_t^{t+s}\stackrel{d}{=}S_{sl}$, where $(S_x)_{x\geq 0}$ is a subordinator with no drift and Lévy measure $\nu$. Adding
that $S_t^{t+s}$ is independent of $(\w{\tau}_{x})_{x\geq0}$ gives  the law of $C^2(t+s)-C^2(t)$. As $(\w{\tau}_{x})_{x\geq0}$  and $S_t^{t+s}$ are
independent of  $(g(u),d(u))_{u \in [0,t]}$, so is  $C^2(t+s)-C^2(t)$. \\ 

These properties ensure that $(g(t),d(t))_{t\in[0,\t{1/m}[}$ is a Markov process. 
\end{proof}

$\newline$
$\bold{} \ \ \ $ To prove Corollary \ref{rate}, we need the following result which uses notation of Theorem \ref{couple}.
\begin{Lem} \label{cvvague} We have the following vague convergence of measure  on $]0,\infty[$ :
$$h^{-1}\P\big(\taup^{(t)}_{S_{hl}}\in \emph{d}x \big) \stackrel{v}{\Longrightarrow } 
l\int_0^{\infty}\nu(\emph{d}z)\P(\taup^{(t)}_z \in \emph{d}x).$$
\end{Lem}
\begin{proof}
Denoting by $\phi$ 
 the Laplace exponent of $(S_x)_{x\geq 0}$, 
$\big(\taup^{(t)}_{S_{xl}}\big)_{x\geq 0}$
is a subordinator of Laplace exponent $l\phi\circ \kappa^{(t)}$ (see (\ref{defexplap})). Moreover for every $\lambda\geq0$,
$\phi(\lambda)=\int_0^{\infty}(1-e^{-\lambda y})\nu(\t{d}y)$, which entails that 
\Bea
\phi\circ \kappa^{(t)}(\lambda)&=& \int_0^{\infty} \big(1-e^{-z\kappa^{(t)}(\lambda)}\big)\nu(\t{d}z)\\
&=&\int_0^{\infty} \E\big(1-e^{-\lambda \taup^{(t)}_z}\big)\nu(\t{d}z)\\
&=&\int_0^{\infty} (1-e^{-\lambda x})\int_0^{\infty}\nu(\t{d}z)\P(\taup^{(t)}_z \in \t{d}x).
\Eea 
Then  $\big(\taup^{(t)}_{S_{xl}}\big)_{x\geq 0}$  is a subordinator with no drift and Lévy measure 
$$l\int_0^{\infty}\nu(\t{d}z)\P(\taup^{(t)}_z \in \t{d}x).$$
Using  Exercise 1 Chapter I in \cite{lev} or $\cite{sub}$ on page 8 completes the proof.
\end{proof}

$\newline$

\begin{proof}[Proof of  Corollary \ref{rate}] We consider first the case when the increment of the left extremity is zero. \\ \\
$\bullet $ Let $c>0$ such that 
$\int_0^{\infty}\nu(\t{d}z)\P(\taup^{(t)}_{z}=c)=0$. Using Theorem $\ref{couple}$ and recalling that
$N_t^{t+h}=N_{]t+t+h]\times\RRR_+^2}=\t{card}\{ i \in \N : T_i \in ]t,t+h] \}$, we have  
\be
\label{tauxcoup}
P(g(t+h)-g(t)=0, \ d(t+h)-d(t) \geq c \ \mid \ l(t)=l)=\P(N_t^{t+h}=0)\P(\taup^{(t)}_{S_{hl}} \geq c).
\ee
Adding that  $\P(N_t^{t+h}=0)\stackrel{h\rightarrow 0}{\longrightarrow }1$ and using Lemma \ref{cvvague} give 
\be
\label{premlim}
h^{-1}P(g(t+h)-g(t)=0, \ d(t+h)-d(t) \geq  c \ \mid \ l(t)=l) \stackrel{h\rightarrow 0}{\longrightarrow } 
l\int_0^{\infty}\nu(\t{d}z)\P(\taup^{(t)}_{z} \geq c).
\ee
$\bullet$ \ Let $a,b>0$  and write
$$P(t,t+h)=\P(g(t)-g(t+h)\geq a, \ d(t+h)-d(t) \geq b \ \mid \ l(t)=l).$$
By Proposition \ref{pppd},  $\big\{(T_i,G_i,D_i) : i \in \N \big\}$ is a PPP on $[0,\t{1/m}]\times \RRR_+^2$ with intensity
d$t\mu^{(t)}(\t{d}y\t{d}x)$. The latter verifies $\P(N_t^{t+h}>1)=o(h)$ ($h\rightarrow 0$), so we have
\be
\label{3}
h^{-1}\P(C^1(t+h)-C^1(t) \in ]-\infty,-a]\times [b,\infty])\stackrel{h\rightarrow 0}{\longrightarrow} \mu^{(t)}([a,\infty[\times [b,\infty]).
\ee
We can prove now that
\be
\label{seclim}
\lim_{h\rightarrow 0} h^{-1}P(t,t+h)= \mu^{(t)} ([a,\infty[\times [b,\infty[).
\ee

- First we give the lowerbound.
$$P(t,t+h)\geq \P(C^1(t+h)-C^1(t) \in ]-\infty,-a]\times [b,\infty] \ \mid \ l(t)=l)$$
Using that $C^1(t+h)-C^1(t)$ is independent of $l(t)$ and  (\ref{3}), we get
\be
\label{1}
\liminf_{h\rightarrow 0} h^{-1}P(t,t+h)\geq  \mu ^{(t)}([a,\infty[\times [b,\infty]).
\ee

-  For the upperbound, observe that
\Bea
P(t,t+h) &\leq& \P(C^1(t+h)-C^1(t) \in ]-\infty,-a]\times [b-\epsilon,\infty] \ \mid \ l(t)=l) \\
&& + \ \ \P(N_t^{t+h}\geq 1, \ C^2(t+h)-C^2(t)\in \{0\}\times[\epsilon,\infty[ \ \mid \ l(t)=l).
\Eea
Using again $C^1(t+h)-C^1(t)$ is independent of $l(t)$ with (\ref{3}) and Theorem $\ref{couple}$ gives
$$
\limsup _{h\rightarrow 0} h^{-1}P(t,t+h) \leq \mu^{(t)}([a,\infty[\times [b-\epsilon,\infty]).
$$
Letting $\epsilon$ tend to $0$ gives the upperbound :
$$
\limsup_{h\rightarrow 0} h^{-1}P(t,t+h) \leq  \mu^{(t)} ([a,\infty[\times [b,\infty[).
$$

The two limits (\ref{premlim}) and (\ref{seclim}) ensure the convergence of measures for sets of the form $\{0\}\times [c,d[$ (with $c>0$)
and $[a,b[\times [c,d[$ (with $a>0$), which completes the proof.
\end{proof}
$\newline$
\section{Evolution  of the right extremity and of the length}
Proposition $\ref{pppd}$, Theorem $\ref{couple}$ and Corollary $\ref{rate}$ give by projection :
\begin{Cor}
\label{d}
$(d(t))_{t\in[0,1/\emph{m}[}$ is a jump process satisfying  \\ \\
(i) $\big\{(T_i,D_i) : i \in \N \big\}$ is a PPP on $[0,\emph{1/m}[\times \RRR^+$ with intensity
$$ \frac{\emph{d}t \int_{z\in [0,\infty]} \emph{d}z \bar{\nu}(z) \P(\taup^{(t)}_z \in \emph{d} x)}{1-\emph{m}t},$$
and $\big\{(T_i,D_i) : i \in \N, \ T_i>t \big\}$ is independent of $(d(u))_{u\in[0,t]}$.  \\ \\
(ii) For all $0\leq t\leq t+s<\emph{1/m}$  : \\
Conditionally on   $l(t)=l$, $d(t+s)-d(t)$ is independent of $(d(u))_{u\in[0,t]}$.    \\
Conditionally  also on $T_i\leq t\leq t+s<T_{i+1}$ for some $i \in \N$ :
$$d(t+s)-d(t)\stackrel{d}{=}\taup^{(t+s)}_{S_{sl}},$$
where  $(S_x)_{x\geq 0}$ is a subordinator with no drift and  Lévy measure $\nu$, that is  independent of $(\taup^{(t+s)}_x)_{x\geq 0}$. \\ \\
The jump rate of $(d(t))_{t\in [0,\emph{1/m}[}$ is given by the following vague convergence of measures on $]0,\infty[$ for $h$ tending to $0$ :
$$\frac{\P(d(t+h)-d(t) \in \emph{d}x \ \mid \ l(t)=l)}{h} \ \stackrel{w}{\Longrightarrow } \
\frac{\int_0^{\infty} \emph{d}z \bar{\nu}(z) \P(\taup^{(t)}_z \in\emph{d}x)}{1-\emph{m}t}+l\int_0^{\infty} \nu(\emph{d}z) \P(\taup^{(t)}_{z}\in \emph{d}x).$$
\end{Cor}
$\newline$
We stress that  $(d(t))_{t \in [0,\t{1/m}[}$ is not a Markov process since the jumps $D_i$ before time $t$ give informations about
$l(t)$ and thus about the future of the process.\\
Note also that we can derive the law of $d(t)$ conditionally on $l(t)$ using Theorem \ref{loi}. More precisely, conditionally
on  $l(t)=l$,
$$\forall d>0 , \qquad  \P(l(t) \in \t{d}l \ \mid \  d(t)=d)=\mathbf{1}_{l\geq d} \frac{\Pi^{(t)}(\t{d}l)}{\bar{\Pi}^{(t)}(d)}.$$
$\newline$

Finally we turn our interest to the process of the length  $(l(t))_{t\in [0,\t{1/m}]}$. Its increments  which are
 due to files arrived
at the left of $g(t)$   which are not stored entirely at the left $g(t)$, are denoted by $L_i$ :
$$L_i:=l(T_i)-l(T_i^-)=G_i+D_i.$$
The other increments of $(l(t))_{t\in [0,\t{1/m}]}$ are due to files which arrive on $\B$.
We can view $(l(t))_{t \in [0,\t{1/m}]}$ as a branching process in continuous time with immigration $L_i$ at time $T_i$
(with no death,  inhomogeneous branching and inhomogeneous immigration) : \\

\begin{Cor}  \label{branch} $(l(t))_{ t \in [0,\emph{1/m}[}$ is an inhomogeneous pure jump Markov process satisfying \\ \\
(i) $\{(T_i,L_i) : i \in \N\}$ is a PPP on $[0,\emph{1/m}[\times \RRR^+$ with intensity
$$\emph{d}t \int_{z\in[0,\infty]}\nu(\emph{d}z)\P(\taup^{(t)}_z\in \emph{d}x)x,$$  and
 $\{(T_i,L_i) : i \in \N, T_i>t\}$ is independent of $(l(s))_{ s \in [0,t]}$ \\ \\
(ii) Conditionally on $T_i\leq t\leq t+s< T_{i+1}$ for some $i\in\N$, $(l(t+u))_{u \in [0,t-s]}$ satisfies the  branching property  :
the law of $(l(t+u))_{u \in [0,t-s]}$  conditioned on $l(t)=x+y$ is equal to the law of the sum of two
independent processes whose laws are   respectively equal
to $(l(t+u))_{u \in [0,t-s]}$  conditioned on $l(t)=x$ and  $(l(t+u))_{u \in [0,t-s]}$  conditioned on $l(t)=y$. \\ \\
The jump rate of $(l(t))_{t \in [0,\emph{1/m}[}$ is given by the following vague convergence of measures on $]0,\infty[$ for $h$ tending to $0$ :
$$\frac{\P(l(t+h)-l(t) \in \emph{d}x \ \mid \ l(t)=l)}{h} \ \stackrel{w}{\Longrightarrow} \ (x+l)\int_0^{\infty}\nu(\emph{d}z)\P(\taup^{(t)}_z\in \emph{d}x).$$
\end{Cor}
$\newline$
\begin{ex}
For the basic example $\nu=\delta_1$, the jump rate of the lenght is equal to
$$\sum_{n=1}^{\infty} \frac{n+l}{n}e^{-tn}\frac{(tn)^{n-1}}{(n-1)!}\delta_n(\t{d}x).$$
This is a consequence of the last displayed limit and (\ref{dirac}).
\end{ex}
$\newline$
\begin{proof}[Proof of  Corollary $\ref{d}$]
Using (\ref{calcint}),  we get :
$$ \int_{z\in[0,\infty]} \P(\taup^{(t)}_z\in \t{d}x)\int_{y\in [0,\infty]}\rho^{(t)}(\t{d}y\t{d}z)= \frac{\int_{z\in[0,\infty]} \t{d}z \bar{\nu}(z) \P(\taup^{(t)}_z \in\t{d}x)}{1-\t{m}t},$$
which gives the intensity of $\big\{(T_i,D_i) : i \in \N \big\}$ by Proposition \ref{pppd}.
\end{proof}

$\newline$

\begin{proof}[Proof of  Corollary $\ref{branch}$]
$(i)$ \ Writing $L_i=G_i+D_i$, Proposition \ref{pppd} entails that $\big\{(T_i,L_i) : i \in \N \big\}$ is a PPP on $[0,\t{1/m}]\times \RRR^{+}$
with intensity $\t{d}t\w{\mu}_t(\t{d}x)$ where
$\w{\mu}_t$ is a measure on $\RRR^+$ defined for a Borel set $A$ of $\RRR^+$ by
$$\w{\mu}_t(A)=\int_{\RRR_+^2}\ind _{\{y+y'\in A\}} \int_0^{\infty} \P(\taup^{(t)}_z\in \t{d}y')\rho^{(t)}(\t{d}y\t{d}z).$$
To determine $\w{\mu}_t$, we compute its  Laplace transform  using Lemma $\ref{calc}$ :

\setlength\arraycolsep{1pt}
\bea
\int_0^{\infty} e^{-\lambda x} \w{\mu}_t(\t{d}x)  &=& \int_{\RRR^{+3}}e^{-\lambda (y+y')}\rho^{(t)}(\t{d}y\t{d}z)\P(\taup^{(t)}_z\in \t{d}y') \nonumber \\
&=&\int_{\RRR^{+3}}e^{-\lambda y'}\P(\taup^{(t)}_z\in \t{d}y')\t{d}z\int_z^{\infty} \nu(\t{d}l) \big[e^{-\lambda y}\P(\taup^{(t)}_{l-z} \in \t{d}y) \nonumber \\
& & \ \ \ \    +\int_0^{y}e^{-\lambda x}\P(\taup^{(t)}_{l-z} \in \t{d}x)(y-x)e^{-\lambda (y-x)}\Pi^{(t)}(\t{d}y-x)\big]  \nonumber \\
&=&\int_0^{\infty}\t{d}ze^{-z\kappa^{(t)}(\lambda)} \int_z^{\infty} \nu(\t{d}l)e^{-(l-z)\kappa^{(t)}(\lambda )}\big[1+\int_0^{\infty}e^{-\lambda u}u\Pi^{(t)}(\t{d}u)\big] \nonumber \\
&=&\int_0^{\infty}\nu(\t{d}l)le^{-l\kappa^{(t)}(\lambda)}[\kappa^{(t)}]'(\lambda )   \nonumber \\
&=& - \frac{\partial  }{\partial  y }\bigg[\int_0^{\infty}\nu(\t{d}l)e^{-l\kappa^{(t)}(y )}\bigg] (\lambda)\nonumber \\ 
&=& - \frac{\partial  }{\partial  y}\bigg[\int_0^{\infty}e^{-y x}\int_0^{\infty}\nu(\t{d}l)\P(\taup^{(t)}_l\in \t{d}x)  \bigg] (\lambda) \nonumber \\
&=& \int_0^{\infty}e^{-\lambda x}x\int_0^{\infty}\nu(\t{d}l)\P(\taup^{(t)}_l\in \t{d}x). \nonumber
\eea
Then $\w{\mu}_t(\t{d}x)=x\int_0^{\infty}\nu(\t{d}z)\P(\taup^{(t)}_z\in \t{d}x)$, which gives the intensity of $\big\{(T_i,L_i) : i \in \N \big\}$. \\ \\
$(ii)$ \ The branching property can be seen as a consequence of the
determination of the jump rate. We give here a more intuitive
approach : We condition by $l(t)=x+y$ and  by $T_i\leq t\leq
t+s<T_{i+1}$ and we make the decomposition effective by splitting
$\B (t)$ in two segments of length $x$ and $y$. First we store the
files $\{i : t_i\in]t,t+s],x_i>d(t)\}$. The free space of the half
line at the right of $\B (t)$  is now the closed range a
subordinator distributed like $(\taup_x^{(t+s)})_{x\geq 0}$ (see
STEP1 in the proof of Corollary \ref{d}).  Then we store
successively the files $\{i : t_i\in]t,t+s],x_i\in[g(t),g(t)+x]\}$
and $\{i : t_i\in]t,t+s],x_i\in]g(t)+x,d(t)]\}$ which induce two
successive increments of the length. The free space at the right
of $0$ after the first storage keeps the same distribution and is
independent of the first increment  by strong regeneration. So the
two increments are independent and distributed respectively like
$l(t+s)-l(t)$ conditioned by $l(t)=x$ and  by $l(t)=y$. This gives
the result since $l(t)$ is Markovian. Formally $l(t+s)-l(t)$ is
equal to   $\taup^{(t+s)}_{S_{s(x+y)}}$ (see proof of  Proposition
$\ref{pppd}$) and
$$\taup^{(t+s)}_{S_{s(x+y)}}= \taup^{(t+s)}_{S_{sx}} \ + \ \taup^{(t+s)}_{S_{s(x+y)}}-\taup^{(t+s)}_{S_{sx}}$$
gives the decomposition expected since $\taup^{(t+s)}_{S_{s(x+y)}}-\taup^{(t+s)}_{S_{sx}}\stackrel{d}{=}
\taup^{(t+s)}_{S_{sy}}$. \\

Using Corollary \ref{rate} and recalling the definition of 
$\w{\mu}_t$ given  at the beginning of the proof ensures that $h^{-1}\P(l(t+h)-l(t) \in \t{d}x  \mid  l(t)=l)$ converges to
$$\w{\mu}_t(\t{d}x)+l\int_0^{\infty}\nu(\t{d}z)\P(\taup^{(t)}_{z}\in \t{d}x).$$
The completes the proof, since $\w{\mu}$ has been determined above.
\end{proof}

$\newline$
\section{Complements}

\subsection{Distribution of $\{(T_i,G_i) : i\in \N\}$ derived from Theorem $\ref{fond}$}

In Section 5, we used the total intensity of the PPP $\{(T_i,G_i) : i\in \N\}$ to prove that the intensity of the PPP
 $\{(T_i,G_i,R_i) : i\in \N\}$ is equal to d$t\rho^{(t)}($d$y$d$z)$ (Theorem \ref{fond}). Here we check that
integrating this intensity with respect to the third coordinate enables us to  recover   the  intensity
of   $\{(T_i,G_i) : i\in \N\}$ given in Theorem \ref{g}.   \\
For that purpose, use Lemma \ref{calc}  to rewrite $ \rho^{(t)}$  as
$$ \rho^{(t)}(\t{d}y\t{d}z)=  \t{d}z\int_0^{\infty} \nu(\t{d}l+z) \big(\P(\taum^{(t)}_{l} \in \t{d}y)+\int_0^{y}\P(\taum^{(t)}_{l} \in \t{d}x)(y-x)\Pi^{(t)}(\t{d}y-x)\big)$$
 and calculate the Laplace transform of   $\int_{z\in[0,\infty]}\rho^{(t)}(\t{d}y\t{d}z)$.
\bea
& & \int_{y\in[0,\infty]}e^{-\lambda y}\int_{z\in [0,\infty]}\rho^{(t)}(\t{d}y\t{d}z) \nonumber \\
&=&\int_0^{\infty } \int_0^{\infty}\t{d}z\nu(\t{d}l+z) \int_0^{\infty}e^{-\lambda y}\big[\P(\taum^{(t)}_{l} \in \emph{d}y)+\int_0^{y}\P(\taum^{(t)}_{l} \in \t{d}x)(y-x)\Pi^{(t)}(\t{d}y-x)\big] \nonumber \\
&=& \int_0^{\infty}\t{d}l \bar{\nu}(l) \big[ e^{-l\kappa(\lambda )}+\int_0^{\infty}\P(\taum^{(t)}_{l} \in \t{d}x)e^{-\lambda x}\int_x^{\infty}e^{-\lambda (y-x)}(y-x)\Pi^{(t)}(\t{d}y-x)\big] \nonumber \\
&=& \int_0^{\infty}\t{d}l \bar{\nu}(l) e^{-l\kappa (\lambda )} [\kappa ^{(t)}]'(\lambda) \nonumber  \\
&=& \int_0^{\infty}\t{d}l \frac{\bar{\nu}(l)}{l} \frac{\partial  }{\partial  \lambda }\E(-e^{-l\kappa^{(t)}(\lambda)})  \nonumber \\
&=& \int_0^{\infty}\t{d}l \frac{\bar{\nu}(l)}{l} \frac{\partial  }{\partial  \lambda }\E(-e^{-\lambda \taum^{(t)}_l}) \nonumber \\
&=& \int_0^{\infty}\t{d}l \frac{\bar{\nu}(l)}{l}\int_0^{\infty}e^{-\lambda y}y\P(\taum_l^{(t)}\in \t{d}y) \nonumber \\
&=& \int_0^{\infty}\t{d}y e^{-\lambda y}\int_0^{\infty}\P(Y^{(t)}_y \in -\t{d}l)\bar{\nu}(l) \ \ \ \t{using} \ (\ref{egal}). \nonumber
\eea
Thus, we conclude with
$$\t{d}t \int_{z\in[0,\infty]}\rho^{(t)}(\t{d}y\t{d}z)= \t{d}t\t{d}x\int_0^{\infty} \P(Y^{(t)}_x \in -\t{d}l)\bar{\nu} (l).$$
$\newline$
\subsection{Direct proof of Corollary $\ref{R}$ using fluctuation theory}
Here we determine the distribution of the remaining data using
fluctuation theory : we get laws at fixed times and  do not need  Theorem \ref{g}, as for the proof of Section 5. \\

 We fix $t$,$h$ and $x$ $\geq0$ . We add the lengths of
files fallen in   $[g(t)-x,g(t)]$ during the time interval $]t,t+h]$. Then
we remove the free space in $[g(t)-x,g(t)]$ at time $t$ which is equal to
$L^{(t)}_x$. The sum of data arrived at the left of $\B (t)$ not
stored at the left of $\B (t)$ between time $t$ and $t+h$ is equal
to the maximum in $x\geq 0$ of this difference. It is also the
quantity of data which has tried to occupy the location $g(t)$
(successfully or not) between time $t$ and $t+h$ :
$Y^{(t+h)}_{g(t)}-I^{(t+h)}_{g(t)}$. So, we have
$\newline$

\begin{Lem}
\label{ppp}
Let $0\leq t<\emph{1/m}$ and $h\geq 0$, then
$$ Y^{(t+h)}_{g(t)}-I^{(t+h)}_{g(t)}\stackrel{} {=}
\emph{sup}\{S_{hx}-L^{(t)}_x, x\geq 0\} \stackrel{} {=}\emph{sup}\{S_{h\taum^{(t)}_x}-x, x\geq 0\} \qquad \t{a.s},$$
where $(S_x)_{x\geq 0}$ is a subordinator with drift $\emph{d}=0$ and Lévy measure $\nu(\t{d}x)$, which is independent of $(L^{(t)}_x)_{x\geq 0}$ and $(\taum^{(t)}_x)_{x\geq 0}$.
\end{Lem}
$\newline$

Denoting  $S^{(t,h)}:=\t{sup}\{S_{h\taum^{(t)}_x}-x, x\geq 0\}$, we have for all $0<a\leq b$,
$$\lim_{h\rightarrow 0} h^{-1} \P(S^{(t,h)}\in [a,b])=\lim_{h\rightarrow 0} h^{-1}\P(\exists  i \in \N : (T_i,R_i) \in ]t,t+h]\times [a,b])$$
and we find the law given in Corollary \ref{R}  :

$\newline$
\begin{Pte} We have the following weak convergence of bounded measures on $]0,\infty[$ when $h$ tends to $0$ :
$$\frac{\P(S^{(t,h)} \in \emph{d}x)}{h} \stackrel{w}{\Longrightarrow }  \frac{ \bar{\nu}(x)\emph{d}x}{1-\emph{m}t}.$$
\end{Pte}
\begin{proof}
$(S_{h\taum^{(t)}_x}-x)_{x\geq 0}$ is a lévy process with negative drift  $-1$, no negative jumps and bounded variation. Its
Laplace exponent is $\kappa^{(t)}\circ (h\phi)-id$, where $\phi$ is the Laplace exponent of $S$ and is defined by
$$\forall \lambda\geq 0, \quad \phi(\lambda)=\int_0^{\infty} (1-e^{-\lambda x})\nu(\t{d}x).$$
Note also that using (\ref{rel}), we have
\bea
\label{rela}
\big[\kappa^{(t)}\circ (h\phi)-id\big]'(0)=[\kappa^{(t)}]'(0).h.\phi'(0)-1=\frac{1}{1-\t{m}t} \t{m} h -1,
\eea
which is negative since $0\leq t+h<1/\t{m}$. Then identity $(14)$ in $\cite{vb}$ or Theorem 5 in \cite{lev}  ensure that $ \forall \lambda >0, \ \forall h \in [0,\t{1/m}-t[$,
$$\E\big( \t{exp} (-\lambda S^{(t,h)})\big)=\bigg(\frac{1}{1-\t{m}t} \t{m} h -1\bigg) \frac{\lambda }{(\kappa^{(t)}\circ (h\phi)-id)(\lambda)} $$
Moreover,
$$\frac{(\kappa^{(t)}\circ (h\phi)-id)(\lambda)}{\lambda }=\frac{\kappa^{(t)}(h\phi(\lambda))}{h\phi(\lambda)}\frac{h\phi(\lambda)}{\lambda}-1=-1+\frac{1}{1-\t{m}t}\frac{h\phi(\lambda)}{\lambda}+\circ _{h\rightarrow 0}(h).$$
So
$$\E\big( \t{exp} (-\lambda S^{(t,h)})\big)=1+\frac{1}{1-\t{m}t}\big(\frac{\phi(\lambda)}{\lambda}-\t{m}\big)h+\circ _{h\rightarrow 0}(h).$$

We can now prove the convergence of $h^{-1}\P(S^{(t,h)}> x)$ when $h$ tends to $0$.
$$\lim_{h\rightarrow 0} \int_0^{\infty}e^{-\lambda x}\frac{\P(S^{(t,h)} > x)}{h}\t{d}x =\lim_{h\rightarrow 0} \frac{1-\E\big( \t{exp} (-\lambda S^{(t,h)})\big)}{h\lambda}
=\frac{1}{1-\t{m}t}\big(\frac{\t{m}}{\lambda} -\frac{\phi(\lambda)}{\lambda^2}\big).$$
Moreover Fubini gives
$$\int_0^{\infty}\t{d}xe^{-\lambda x}\int_x^{\infty} \bar \nu (a)\t{d}a =
\int_0^{\infty}\nu(\t{d}y)\int_0^y \t{d}a \frac{1-e^{-\lambda a}}{\lambda}
=\frac{\t{m}}{\lambda} -\frac{\phi(\lambda)}{\lambda^2}. $$
Then for every $\lambda>0$,
$$\lim_{h\rightarrow 0} \int_0^{\infty}e^{-\lambda x}\frac{\P(S^{(t,h)} > x)}{h}\t{d}x =
\int_0^{\infty}e^{-\lambda x}\frac{\int_x^{\infty} \bar \nu (a)\t{d}a}{1-\t{m}t}\t{d}x,
$$
which proves the convergence of
$\P(S^{(t,h)} \in \t{d}x)/h$ to $\bar{\nu}(x)\t{d}x/(1-\t{m}t)$. Indeed, introduce
the measures  $\mu_h(\t{d}x)$ and $\mu(\t{d}x)$ on $\RRR^+$ whose tails are given by
$$\mu_h(]x,\infty])=e^{-x}\P(S^{(t,h)} > x)/h, \qquad
\mu (]x,\infty])=e^{-x}\int_x^{\infty} \bar{\nu}(a)\t{d}a/(1-\t{m}t).$$
The last displayed limit entails the weak convergence of $\mu_h(\t{d}x)$ to  $\mu(\t{d}x)$ when $h$ tends to $0$,
by convergence of  Laplace transforms. As $\mu$ is non atomic, for every $x\geq 0$,
$\mu_h(]x,\infty])$ tends to $\mu (]x,\infty])$, which proves that
$\P(S^{(t,h)} > x)/h$ tends to $\int_x^{\infty} \bar{\nu}(a)\t{d}a/(1-\t{m}t).$
\end{proof}

\begin{Rque}
Denote $\gamma^{(t,h)}$ the a.s instant at which the supremum
$S^{(t,h)}$ is reached. To obtain the distribution of
$\{(T_i,G_i,R_i)\  : \ i \in \N\}$ by this way, we need to know
the joint law of $(S^{(t,h)},\taum^{(t)}_{\gamma^{(t,h)}})$ which
we cannot derive directly from fluctuation theory.
\end{Rque}

$\newline$

\end{document}